\documentclass{article}

\usepackage[final]{neurips_2019}

\usepackage[utf8]{inputenc} 
\usepackage[T1]{fontenc}    
\usepackage{hyperref}       
\usepackage{url}            
\usepackage{booktabs}       
\usepackage{amsfonts}       
\usepackage{nicefrac}       
\usepackage{microtype}      

\usepackage{algorithm}
\usepackage{algorithmic}

\usepackage{epsfig}
\usepackage{amssymb}
\usepackage{amsmath}
\usepackage{amsthm}
\usepackage{amsfonts}
\usepackage{bbding}
\usepackage{array}
\usepackage{caption,tabularx,booktabs}
\usepackage{arydshln}
\usepackage{xargs}                      
\usepackage[pdftex,dvipsnames]{xcolor}  

\usepackage[colorinlistoftodos,prependcaption,textsize=tiny]{todonotes}
\newcommandx{\unsure}[2][1=]{\todo[inline,linecolor=red,backgroundcolor=red!25,bordercolor=red,#1]{#2}}
\newcommandx{\change}[2][1=]{\todo[linecolor=blue,backgroundcolor=blue!25,bordercolor=blue,#1]{#2}}
\newcommandx{\info}[2][1=]{\todo[linecolor=OliveGreen,inline,backgroundcolor=OliveGreen!25,bordercolor=OliveGreen,#1]{#2}}
\newcommandx{\improvement}[2][1=]{\todo[linecolor=Plum,inline,backgroundcolor=Plum!25,bordercolor=Plum,#1]{#2}}

\usepackage[font=small,skip=5pt]{caption}

\newtheorem{thm}{Theorem}

\newtheorem{ass}{Assumption}

\newcolumntype{C}[1]{>{\centering\let\newline\\\arraybackslash\hspace{0pt}}m{#1}}

\newcommand\tagthis{\addtocounter{equation}{1}\tag{\theequation}}

\DeclareMathOperator{\Ocal}{\mathcal{O}}

\makeatletter
\newcounter{subthm} 
\let\savedc@thm\c@hyp

\newcommand{\normhyp}{%
  \let\c@hyp\savedc@hyp 
  \renewcommand\thehyp{\arabic{hyp}}%
} 
\makeatother

\makeatletter
\newcounter{subass} 
\let\savedc@ass\c@hyp

\makeatother

\title{Tight Dimension Independent Lower Bound on the Expected Convergence Rate for Diminishing Step Sizes in SGD 
}

\author{%
  Phuong Ha Nguyen \\
  Electrical and Computer Engineering \\
            University of Connecticut, USA \\
  \texttt{phuongha.ntu@gmail.com} \\
  \And
  Lam M. Nguyen \\
  IBM Research, Thomas J. Watson Research Center \\
            Yorktown Heights, USA \\
  \texttt{LamNguyen.MLTD@ibm.com} \\
  \And
  Marten van Dijk \\
  Electrical and Computer Engineering \\
            University of Connecticut, USA \\
  \texttt{marten.van$\_$dijk@uconn.edu} \\  
}

\begin{document}

\maketitle

\begin{abstract}
We study the convergence of Stochastic Gradient Descent (SGD) for strongly convex  objective functions. We prove for all $t$ a lower bound on the expected convergence rate after the $t$-th SGD iteration; the lower bound is over all possible sequences of diminishing step sizes. It implies that recently proposed sequences of step sizes at ICML 2018 and ICML 2019 are {\em universally} close to optimal in that the expected convergence rate after {\em each} iteration is within a factor $32$ of our lower bound. This factor is independent of dimension $d$. We offer a framework for comparing with lower bounds in state-of-the-art literature and when applied to SGD for strongly convex objective functions our lower bound is a significant factor $775\cdot d$ larger compared to existing work. 


\end{abstract}

\section{Introduction}

We are interested in solving the following stochastic optimization problem
\begin{align*}
\min_{w \in \mathbb{R}^d} \left\{ F(w) = \mathbb{E} [ f(w;\xi) ] \right\}, \tagthis \label{main_prob_expected_risk}  
\end{align*}
where $\xi$ is a random variable obeying some distribution $g(\xi)$. 
In the case of empirical risk minimization with a training set $\{(x_i,y_i)\}_{i=1}^n$, $\xi_i$ is a random variable that is defined 
 by a single random sample  $(x,y)$  pulled uniformly from the training set. Then,  by defining  $f_i(w) := f(w;\xi_{i})$, empirical risk minimization reduces to
\begin{gather}\label{main_prob}
\min_{w \in \mathbb{R}^d} \left\{ F(w) = \frac{1}{n} \sum_{i=1}^n f_i(w) \right\}.  
\end{gather}

Problems of this type arise frequently in supervised learning applications \cite{ESL}. The classic first-order methods to solve  problem (\ref{main_prob}) are gradient descent (GD) \cite{Nocedal2006NO}  and stochastic gradient descent (SGD)\footnote{We notice that even though stochastic gradient is referred to as SG in literature, the term stochastic gradient descent (SGD) has been widely used in many important works of large-scale learning.} \cite{RM1951} algorithms. GD is a standard deterministic gradient method, which updates iterates along the negative full gradient with  learning rate $\eta_t$ as follows
\begin{gather*}
w_{t+1} = w_{t} - \eta_t \nabla F(w_{t}) = w_{t} - \frac{\eta_t}{n} \sum_{i=1}^n \nabla f_i(w_{t}) \ , \ t \geq 0.  
\end{gather*}
We can choose $\eta_t = \eta = \Ocal(1/L)$ and achieve a linear convergence rate for the strongly convex case \cite{nesterov2004}. 
The upper bound of the convergence rate of GD and SGD has been studied in~ \cite{Bertsekas1999,BoydVandenberghe2004,nesterov2004,schmidt2013fast,Nguyen2018,Nguyen2018_NewAspectsSGD,gower2019sgd}. 

The disadvantage of GD is that it requires evaluation of $n$ derivatives at each step, which is very expensive and therefore avoided in large-scale optimization. To reduce the computational cost for solving \eqref{main_prob}, a class of variance reduction methods \cite{SAG,SAGA,SVRG,Nguyen2017_sarah} has been proposed. The difference between GD and variance reduction methods is that GD needs to compute the full gradient at each step, while the variance reduction methods will compute the full gradient after a certain number of steps. In this way, variance reduction methods have less computational cost compared to GD. 
To avoid evaluating the full gradient at all, SGD generates an unbiased random variable $\xi_t$ satisfying
$$\mathbb{E}_{\xi_t}[\nabla f(w_t;\xi_t)]=\nabla F(w_t),$$
and then evaluates gradient $\nabla f(w_t;\xi_t)$ for $\xi_t$ drawn from  a distribution $g(\xi)$. After this, $w_t$ is updated as follows
\begin{gather*}
w_{t+1} = w_{t} - \eta_t \nabla f(w_{t};\xi_t). \tagthis \label{sgd_step}
\end{gather*}

We focus on the general problem \eqref{main_prob_expected_risk} where $F$ is strongly convex. Since $F$ is strongly convex, a unique optimal solution of \eqref{main_prob_expected_risk} exists and throughout the paper we denote this optimal solution by $w_{*}$ and are interested in studying the expected convergence rate 
$$Y_t= \mathbb{E}[ \|w_t-w_*\|^2 ].$$ 

Algorithm \ref{sgd_algorithm} provides a detailed description of SGD.
Obviously, the computational cost of a single iteration in SGD is $n$ times cheaper than that of a single iteration in GD. However, as has been shown in literature we need to choose $\eta_t = \Ocal(1/t)$ and the expected convergence rate of SGD is slowed down to $\Ocal(1/t)$ \cite{bottou2016optimization}, which is a sublinear convergence rate. 

\begin{algorithm}[h]
   \caption{Stochastic Gradient Descent (SGD) Method}
   \label{sgd_algorithm}
\begin{algorithmic}
   \STATE {\bfseries Initialize:} $w_0$
   \STATE {\bfseries Iterate:}
   \FOR{$t=0,1,\dots$}
  \STATE Choose a step size (i.e., learning rate) $\eta_t>0$. 
  \STATE Generate a random variable $\xi_t$ with probability density $g(\xi_t)$.
  \STATE Compute a stochastic gradient $\nabla f(w_{t};\xi_{t}).$
   \STATE Update the new iterate $w_{t+1} = w_{t} - \eta_t \nabla f(w_{t};\xi_{t})$.
   \ENDFOR
\end{algorithmic}
\end{algorithm}

\vspace{2mm}

\noindent
{\bf Problem Statement and Contributions:} 
We seek to find
 a tight lower bound on the expected convergence rate
 $Y_t$ with the purpose of showing that the stepsize sequences of \cite{Nguyen2018} and \cite{gower2019sgd} 
 for classical SGD is optimal for $\mu$-strongly convex and $L$-smooth respectively expected $L$-smooth objective functions within a {\em small dimension independent constant factor}. This is important because of the following reasons: 
 \begin{enumerate}
 \item
The lower bound tells us that a sequence of stepsizes as a function of only $\mu$ and $L$ cannot beat an expected convergence rate of $O(1/t)$ -- this is known general knowledge and was already proven in~\cite{AgarwalBartlettRavikumarEtAl}, where a {\em dimension dependent} lower bound for a larger class of algorithms that includes SGD was proven.  
For the class of SGD with diminishing stepsizes as a  function of only global parameters $\mu$ and $L$ we show a {\em dimension independent} lower bound which is a factor $775\cdot d$ larger.
\item  We now understand 
into what extent the sequence of stepsizes of \cite{Nguyen2018} and \cite{gower2019sgd} are 
optimal in that it leads to minimal expected convergence rates $Y_t$ for {\em all} $t$:  For each $t$ we will show a {\em dimension independent} lower bound   on $Y_t$ over {\em all possible} stepsize sequences. This
includes the {\em best possible} stepsize sequence which minimizes $Y_t$ for a {\em given} $t$. Our lower bound achieves the upper bound on $Y_t$ for the stepsize sequences of \cite{Nguyen2018} and \cite{gower2019sgd} 
within a factor 32 for {\em all} $t$. This implies that these stepsize sequences 
universally minimizes each $Y_t$ within factor 32. 
\item As a consequence, in order to attain a better expected convergence rate, we need to {\em either} assume more specific knowledge about the objective function $F$ so that we can construct a better stepsize sequence for SGD based on this additional knowledge {\em or} we need to step away from SGD and use a different kind of algorithm. For example, the larger class of algorithms in~\cite{AgarwalBartlettRavikumarEtAl} may contain a non-SGD algorithm which may get close to the lower bound proved in~\cite{AgarwalBartlettRavikumarEtAl} which is a factor $775\cdot d$ smaller. Since the larger class of algorithms in ~\cite{AgarwalBartlettRavikumarEtAl}  contains algorithms such as Adam~\cite{kingma2014adam}, AdaGrad~\cite{duchi2011adaptive}, SGD-Momentum \cite{sutskever2013importance}, RMSProp \cite{zou2018sufficient}
we now know that these practical algorithms will at most improve a factor $32\cdot 775\cdot d$ over SGD for strongly convex optimization -- this can be significant as this can lead to orders of magnitude less gradient computations. 
We are the first to make such quantification.
\end{enumerate}

\vspace{2mm}

\noindent
{\bf Outline:}
Section \ref{sec:related} discusses background: First, we discuss the recurrence on $Y_t$ used in \cite{Nguyen2018} for proving their upper bound on $Y_t$ -- this recurrence plays a central role in proving our lower bound. We discuss the upper bounds of both \cite{Nguyen2018} and \cite{gower2019sgd} -- the latter holding for a larger class of algorithms. Second, we explain the lower bound of ~\cite{AgarwalBartlettRavikumarEtAl} in detail in order to be able to properly compare with our lower bound. Section \ref{sec:frame} introduces a framework for comparing bounds and explains the consequences of our lower bound in detail. Section~\ref{sec:tight1} describes a class of strongly convex and smooth objective functions which is used to derive our lower bound.
We also verify our theory by experiments in the supplementary material.
Section~\ref{sec:conclusion} concludes the paper.

\section{Background} \label{sec:related}

We explain the upper bound of \cite{Nguyen2018,gower2019sgd}, and lower bound of \cite{AgarwalBartlettRavikumarEtAl} respectively.

\subsection{Upper Bound for Strongly Convex and Smooth Objective Functions}

The starting point for analysis is the recurrence first introduced in \cite{Nguyen2018,Leblond2018} 
\begin{align*}
\mathbb{E}[\| w_{t+1} - w_{*} \|^2 ]  \leq (1-\mu \eta_t) \mathbb{E}[ \| w_{t} - w_{*} \|^2 ] + \eta_t^2 N, \tagthis \label{main_ineq_sgd_new02}
\end{align*}
where  
$$N = 2 \mathbb{E}[ \|\nabla f(w_{*}; \xi)\|^2 ]$$
and $\eta_t$ is upper bounded by $\frac{1}{2L}$; the recurrence has been shown to hold, see \cite{Nguyen2018,Leblond2018}, if we assume
\begin{enumerate}
    \item $F(.)$ is $\mu$-strongly convex, \item $f(w;\xi)$ is $L$-smooth, 
    \item $f(w;\xi)$ is convex, and 
    \item $N$ is finite; 
\end{enumerate}

we detail these assumptions below:

\begin{ass}[$\mu$-strongly convex]
\label{ass_stronglyconvex}
The objective function $F: \mathbb{R}^d \to \mathbb{R}$ is  $\mu$-strongly convex, i.e., there exists a constant $\mu > 0$ such that $\forall w,w' \in \mathbb{R}^d$, 
\begin{gather*}
F(w) - F(w') \geq \langle \nabla F(w'), (w - w') \rangle + \frac{\mu}{2}\|w - w'\|^2. \tagthis\label{eq:stronglyconvex_00}
\end{gather*}
\end{ass}

\begin{ass}[$L$-smooth]
\label{ass_smooth}
$f(w;\xi)$ is $L$-smooth for every realization of $\xi$, i.e., there exists a constant $L > 0$ such that, $\forall w,w' \in \mathbb{R}^d$, 
\begin{align*}
\| \nabla f(w;\xi) - \nabla f(w';\xi) \| \leq L \| w - w' \|. \tagthis\label{eq:Lsmooth_basic}
\end{align*} 
\end{ass}

Assumption \ref{ass_smooth} implies that $F$ is also $L$-smooth. 


\begin{ass}\label{ass_convex}
$f(w;\xi)$ is convex for every realization of $\xi$, i.e., $\forall w,w' \in \mathbb{R}^d$, 
\begin{gather*}
f(w;\xi)  - f(w';\xi) \geq \langle \nabla f(w';\xi),(w - w') \rangle.
\end{gather*}
\end{ass}

\begin{ass}\label{finiteN}
$N = 2 \mathbb{E}[ \|\nabla f(w_{*}; \xi)\|^2 ]$ is finite.
\end{ass}

We denote the set of strongly convex objective functions by ${\cal F}_{str}$ and denote the subset of ${\cal F}_{str}$ satisfying Assumptions \ref{ass_stronglyconvex}, \ref{ass_smooth}, \ref{ass_convex}, and \ref{finiteN} by ${\cal F}_{sm}$.
   
We notice that the earlier established recurrence in \cite{Moulines2011} under the same set of assumptions 
\begin{align*}
\mathbb{E}[\| w_{t+1} - w_{*} \|^2 ]  &\leq (1-2\mu \eta_t + 2L^2\eta_t^2) \mathbb{E}[ \| w_{t} - w_{*} \|^2 ] + \eta_t^2 N
\end{align*}
is similar, but worse than (\ref{main_ineq_sgd_new02}) as it only holds for $\eta_t < \frac{\mu}{L^2}$ where  (\ref{main_ineq_sgd_new02}) holds for $\eta_t\leq \frac{1}{2L}$.  Only for step sizes $\eta_t < \frac{\mu}{2L^2}$  the above recurrence provides a better bound than (\ref{main_ineq_sgd_new02}), i.e., $1-2\mu \eta_t + 2L^2\eta_t^2 \leq 1-\mu \eta_t $. 
In practical settings such as logistic regression  $\mu=\Ocal(1/n)$, $L=\Ocal(1)$, and $t=\Ocal(n)$ (i.e. $t$ is at most a relatively small constant number of epochs, where a single epoch represents $n$ iterations resembling the complexity of a single GD computation). 
See (\ref{step}) below, for this parameter setting the optimally chosen step sizes are $\gg \frac{\mu}{L^2}$. This is the reason we focus in this paper  on analyzing recurrence  (\ref{main_ineq_sgd_new02}) in order to prove our lower bound:
For $\eta_t \leq \frac{1}{2L}$,
\begin{equation}  \hspace{3mm}
 Y_{t+1} \leq (1 - \mu \eta_t)Y_t + \eta_t^2 N, \label{EqY}
 \end{equation}
where $Y_t= \mathbb{E}[ \|w_t-w_*\|^2 ]$.
 




Based on the above assumptions (\textit{without} the so-called bounded gradient assumption) and knowledge of only $\mu$ and $L$ a sequence of step sizes $\eta_t$ can be constructed such that $Y_t$ is smaller than $\Ocal(1/t)$~\cite{Nguyen2018}; more explicitly, 
for the sequence of step sizes 
\begin{equation} \eta_t=\frac{2}{\mu t+4L} \label{step} \end{equation}
we have for all objective functions in ${\cal F}_{sm}$ the upper bound
\begin{equation} Y_t \leq \frac{16N}{\mu}\frac{1}{\mu(t-T')+4L} = \frac{16N}{\mu^2 t}(1+{\cal O}(1/t)),\label{Yt} \end{equation}
where 
$$t\geq T'=\frac{4L}{\mu}\max\{\frac{L\mu Y_0}{N},1\}-\frac{4L}{\mu}.$$ 

We notice that ~\cite{gower2019sgd}  studies the larger class, which we denote ${\cal F}_{esm}$, which is defined as ${\cal F}_{sm}$ where expected smoothness is assumed in stead of smoothness and convexity of component functions. We rephrase their assumption for classical SGD as studied in this paper.\footnote{This means that distribution ${\cal D}$ in \cite{gower2019sgd} must be over unit vectors $v\in [0,\infty)^n$, where $n$ is the number of component functions, i.e., $n$ possible values for $\xi$. Arbitrary distributions ${\cal D}$ correspond to SGD with mini-batches where each component function indexed by $\xi$ is weighted with $v_\xi$.}

\begin{ass} \label{ass-smexp} ($L$-smooth in expectation) The objective function $F: \mathbb{R}^d \to \mathbb{R}$ is  $L$-smooth in expectation if there exists a constant $L > 0$ such that, $\forall w \in \mathbb{R}^d$, 
\begin{align*}
\mathbb{E}[ \| \nabla f(w;\xi) - \nabla f(w_*;\xi) \|^2 ] \leq 2L \| F(w) - F(w_*) \|. \tagthis\label{eq:Lsmoothexp}
\end{align*} 
\end{ass}

The results in \cite{gower2019sgd} assume the above assumption for empirical risk minimization (\ref{main_prob}). $L$-smoothness, see \cite{nesterov2004}, implies Lipschitz continuity (i.e., $\forall w, w' \in \mathbb{R}^d$, 
$$f(w,\xi) \leq f(w',\xi) + \langle \nabla f(w',\xi),(w-w') \rangle + \frac{L}{2}\|w-w'\|^2$$) and together with Proposition A.1 in \cite{gower2019sgd} this implies $L$-smooth in expectation. This shows that ${\cal F}_{esm}$ defined by Assumptions \ref{ass_stronglyconvex}, \ref{finiteN}, and \ref{ass-smexp} is indeed a superset of ${\cal F}_{sm}$.

The step sizes (\ref{step}) from~\cite{Nguyen2018} for ${\cal F}_{sm}\subseteq {\cal F}_{esm}$ and
\begin{equation}
\eta_t = \frac{2t+1}{(t+1)^2\mu} \mbox{ for } t> \frac{4L}{\mu} \mbox{ and } \eta_t=\frac{1}{2L} \mbox{ for } t\leq \frac{4L}{\mu} \label{stepe}
\end{equation}
developed for ${\cal F}_{esm}$ in~\cite{gower2019sgd} and ~\cite{Nguyen2018} are equivalent in that they are both $\approx \frac{2}{\mu t}$ for $t$ large enough. Both step size sequences give exactly the same asymptotic upper bound (\ref{Yt}) on $Y_t$ (in our notation).


In \cite{RM1951}, the authors proved the convergence of SGD for the step size sequence $\{\eta_t\}$ satisfying conditions
$
\sum_{t=0}^{\infty} \eta_t = \infty \ \text{and} \ \sum_{t=0}^{\infty} \eta_t^2 < \infty$.
In \cite{Moulines2011}, the authors studied the expected convergence rates for another class of step sizes of $\mathcal{O}(1/t^p)$ where $0< p\leq 1$.
However, the authors of both \cite{RM1951} and \cite{Moulines2011} do {\em not} discuss about the optimal step sizes among all proposed step sizes which is what is done in this paper.

\subsection{Lower Bound for First Order Stochastic Oracles}




The authors of~\cite{NemirovskyYudin1983} proposed the first formal study on lower bounding the expected convergence rate for a large class of algorithms which includes SGD. The authors of~\cite{AgarwalBartlettRavikumarEtAl} and~\cite{RaginskyRakhlin2011} independently studied this lower bound using information theory and were able to improve it.  


The derivation in~\cite{AgarwalBartlettRavikumarEtAl} is for algorithms including SGD where the sequence of stepsizes is a-priori fixed based on global information regarding assumed stochastic parameters concerning the objective function $F$. 
%
Their proof uses the following set of assumptions: First, 
The assumption of a strongly convex objective function, i.e., Assumption~\ref{ass_stronglyconvex} (see Definition 3 in~\cite{AgarwalBartlettRavikumarEtAl}). Second, the objective function is convex Lipschitz:

\begin{ass} \label{ass:lip} (convex Lipschitz) The objective function $F$ is a convex Lipschitz function, i.e., there exists a bounded convex set $\mathcal{S} \subset \mathbb{R}^d$ and a positive number $K$ such that $\forall w, w' \in \mathcal{S} \subset \mathbb{R}^d$
\begin{equation*}
\|  F(w) - F(w') \| \leq K \|w-w' \|.
\end{equation*} 
\end{ass}

We notice that this assumption implies the assumption on bounded gradients as stated here (and explicitly mentioned in Definition 1 in \cite{AgarwalBartlettRavikumarEtAl}): 
There exists a bounded convex set $\mathcal{S} \subset \mathbb{R}^d$ and a positive number $\sigma$ such that 
\begin{equation}
\mathbb{E}[\|\nabla f(w;\xi) \|^2]\leq \sigma^2 \label{sigma}
\end{equation}
for all $w\in \mathcal{S} \subset \mathbb{R}^d$. This is not the same as the bounded gradient assumption where $\mathcal{S} = \mathbb{R}^d$ is unbounded.\footnote{The bounded gradient assumption, where $\mathcal{S}$ is unbounded, is in conflict with assuming strong convexity as explained in \cite{Nguyen2018}.} Clearly, for $w_*$, (\ref{sigma}) implies a finite $N\leq 2\sigma^2$.

We define ${\cal F}_{lip}$ as the set of strongly convex objective functions that satisfy Assumption \ref{ass:lip}. Classes ${\cal F}_{esm}$ and ${\cal F}_{lip}$ are both subsets of ${\cal F}_{str}$ and differ (are not subclasses of each other) in that they assume expected smoothness and convex Lipschitz respectively.

To prove a lower bound of $Y_t$ for ${\cal F}_{lip}$, the authors constructed a class of objective functions $\subseteq {\cal F}_{lip}$ and showed  a lower bound of $Y_t$ for this class; in terms of the notation used in this paper, 
\begin{equation}
 \frac{\log(2/\sqrt{e})}{432 \cdot d}\frac{N}{\mu^2 t}.\label{eqLB}   
\end{equation}

The authors of~\cite{AgarwalBartlettRavikumarEtAl}  prove lower bound (\ref{eqLB}) for the class ${\cal A}_{stoch}$ of {\em stochastic first order algorithms} that can be understood as operating based on information provided by a stochastic  first-order oracle, i.e., any algorithm which bases its computation in the $t$-th iteration on $\mu$, $K$ or $L$, $d$, and access to an oracle that provides $f(w_t;\xi_t)$ and $\nabla f(w_t;\xi_t)$. This class includes ${\cal A}_{SGD}$ defined as SGD with some sequence of diminishing step sizes as a function of global parameters such as $\mu$ and $L$ or $\mu$ and $K$, see  Algorithm \ref{sgd_algorithm}. We notice that ${\cal A}_{stoch}$ also includes practical algorithms such as Adam~\cite{kingma2014adam}, 
etc.
%
We revisit their derivation in the supplementary material
where we show\footnote{We also discuss the underlying assumption of convex Lipschitz and show that in order for the analysis in \cite{AgarwalBartlettRavikumarEtAl} to follow through one -- likely tedious but believable -- statement still needs a formal proof.} how their lower bound transforms into (\ref{eqLB}). Notice that their lower bound depends on dimension $d$.




\section{Framework for Upper and Lower Bounds} \label{sec:frame}

Let $par(F)$ denote the concrete values of the global parameters of an objective function $F$ such as the values for $\mu$ and $L$ corresponding to objective functions $F$ in ${\cal F}_{sm}$ and ${\cal F}_{esm}$ or $\mu$ and $K$ corresponding to objective functions $F$ in ${\cal F}_{lip}$. When defining a class ${\cal F}$ of objective functions,  we also need to explain how ${\cal F}$ defines a corresponding $par(.)$ function. We will use the notation ${\cal F}[p]$ to stand for the subclass $\{ F \in {\cal F} \ : \ p=par(F)\}\subseteq {\cal F}$, i.e., the subclass of objective functions of ${\cal F}$ with the same parameters $p$. We assume that parameters of a class are included in the parameters of a smaller subclass: For example, ${\cal F}_{sm}$ is a subset of the class of strongly convex objective functions ${\cal F}_{str}$ with only global parameter $\mu$. This means that for concrete values $\mu$ and $L$ we have ${\cal F}_{sm}[\mu,L]\subseteq {\cal F}_{str}[\mu]$.


For a given objective function $F$, we are interested in the best possible expected convergence rate after the $t$-th iteration among all possible algorithms $A$ in a larger class of algorithms ${\cal A}$. Here, we assume that ${\cal A}$ is a subclass of the larger class ${\cal A}_{stoch,{\cal U}}$ of stochastic first order algorithms where the computation in the $t$-th iteration not only has access to $par(F)$ and access to an oracle  that provides $f(w_t;\xi_t)$ and $\nabla f(w_t;\xi_t)$ but also access to possibly another oracle ${\cal U}$ providing even more information. Notice that ${\cal A}\subseteq {\cal A}_{stoch}\subseteq {\cal A}_{stoch,{\cal U}}$ for any oracle ${\cal U}$. With respect to the expected convergence rate, we want to know which algorithm $A$ in ${\cal A}$ minimizes $Y_t$ the most. Notice that for different $t$ this may be a different algorithm $A$. We define for $F\in {\cal F}$ (with associated $par(.)$)
$$\gamma^F_t({\cal A})= \underset{A \in {\cal A}}{\mathrm{inf}} \ Y_t(F,A),$$ where $Y_t$ is explicitly shown as a function of the objective function $F$ and choice of algorithm $A$.

Among the objective functions $F \in {\cal F}$ with same global parameters $p=par(F)$ (i.e., $F\in {\cal F}[p]$), we consider the objective function $F$ which has the \textit{worst} expected convergence rate at the $t$-th iteration. This is of interest to us because algorithms $A$ only have access to $p=par(F)$ as the sole information about objective function $F$, hence, if we prove an upper bound on the expected convergence rate for algorithm $A$, then this upper bound must hold for all $F\in {\cal F}$ with the same parameters $p=par(F)$. In other words such an upper bound must be at least

\begin{align*}
 \gamma_t({\cal F}[p], {\cal A})&=\underset{F \in \mathcal{F}[p]}{\mathrm{sup}} \ \gamma^F_t({\cal A})  =
 \underset{F \in \mathcal{F}[p]}{\mathrm{sup}} \ \ \underset{A\in {\cal A}}{\mathrm{inf}} \ Y_t(F,A).   
\end{align*}

So, any lower bound on $\gamma_t({\cal F}[p], {\cal A})$ gives us a lower bound on the best possible upper bound on $Y_t$ that can be achieved. Such a lower bound tells us into what extent the expected convergence rate $Y_t$ cannot be improved. 

The lower bound (\ref{eqLB}) and upper bound (\ref{Yt}) are not only a function of $\mu$ in $p=par(F)$ but also a function of $N$  which is outside $p=par(F)$ for $F\in {\cal F}_{lip}$ or $F\in {\cal F}_{esm}$. We are really interested in such more fine-grained bounds that are a function of $N$. For this reason we need to consider the subclass of objective functions $F$ in ${\cal F}[p]$ that all have the same $N$. We implicitly understand that $N$ is an auxiliary parameter of an objective function $F$ and we denote this as a function of $F$ as $N(F)$. We define ${\cal F}^a[p]= \{ F\in {\cal F}[p] \ : \ a = aux(F)\}$ where $aux(.)$ represents for example $N(.)$. This leads to notation like ${\cal F}^{N}_{lip}[\mu, K,d]$. Notice that $p=par(F)$ can be used by an algorithm $A\in {\cal A}$ while $a=aux(F)$ is not available to $A$ through $p=par(F)$ (but may be available through access to an oracle).

If we find a tight lower bound with upper bound up to a constant factor, as in this paper, then we know that the algorithm that achieves the upper bound is close to optimal in that the expected convergence rate cannot be further minimized/improved in a significant way. In practice we are only interested in upper bounds on $Y_t$ that can be realized by the {\em same} algorithm $A$ (if not, then we need to know a-priori the exact number of iterations $t$ we want to run an algorithm and then choose the best one for that $t$). In this paper we consider the algorithm $A$  for $F$ in ${\cal F}_{sm}$ resp. ${\cal F}_{esm}$ defined as SGD with diminishing step sizes (\ref{step}) resp. (\ref{stepe}) as a function of $par(F)=(\mu, L)$ giving upper bound (\ref{Yt}) on expected convergence rate $Y_t(F,A)$. We show that $A$ is close to optimal. 

Given the above definitions we have
\begin{equation} \gamma_t({\cal F}[p], {\cal A}) \leq \gamma_t({\cal F}'[p'], {\cal A}') 
\label{gammaineq} \end{equation}
for ${\cal F}[p]\subseteq {\cal F}'[p']$ and ${\cal A}'\subseteq {\cal A}$, i.e., the worst objective function in a larger class of objective functions is worse than the worst objective function in a smaller class of objective functions (see the supremum used in defining $\gamma_t$) and the best algorithm from a larger class of algorithms is better than the best algorithm from a smaller class of algorithms (see the infinum used in defining $\gamma_t$). 
This implies 
\begin{eqnarray} \gamma_t({\cal F}^N_{lip}[\mu,K,d], {\cal A}_{stoch}) &\leq  & \gamma_t({\cal F}^N_{str}[\mu], {\cal A}_{SGD}), \label{FAstoch} \\
\gamma_t({\cal F}^N_{sm}[\mu,L], {\cal A}_{ExtSGD}) &\leq & \gamma_t({\cal F}^N_{esm}[\mu,L], {\cal A}_{SGD}) \leq \gamma_t({\cal F}^N_{str}[\mu], {\cal A}_{SGD}), \label{FAExtSGD}
\end{eqnarray}
where ${\cal A}_{SGD}\subseteq {\cal A}_{ExtSGD}$ is defined as follows:

In our framework we introduce {\em extended SGD} as the class ${\cal A}_{ExtSGD}$ of SGD algorithms where the stepsize in the $t$-th iteration can be computed based on global parameters $\mu$, $L$, and access to an oracle ${\cal U}$ that provides additional information $N$, $\nabla F(w_t)$, and $Y_t$. This class also includes SGD with diminishing stepsizes as defined in Algorithm  \ref{sgd_algorithm}, i.e., ${\cal A}_{SGD}\subseteq {\cal A}_{ExtSGD}$. The reason for introducing the larger class ${\cal A}_{ExtSGD}$ is not because it contains practical algorithms different than SGD, on the contrary. The only reason is that it allows us to define {\em one single algorithm} $A\in {\cal A}_{ExtSGD}$ which realizes $\gamma_t^F({\cal A}_{ExtSGD})$ {\em for all} $t$ for all $F$ in a to be constructed subclass ${\cal F}\subseteq {\cal F}_{sm}$ -- the topic of the next section. This property allows a rather straightforward calculus based proof without needing to use more advanced concepts from information and probability theory as required in the proof of \cite{AgarwalBartlettRavikumarEtAl}. Looking ahead, we will prove in Theorem \ref{theorem:example_tightness}
\begin{equation}
    \frac{1}{2}\frac{N}{\mu^2 t}(1- {\cal O}((\ln t)/t)) \leq \gamma_t({\cal F}^N_{sm}[\mu,L], {\cal A}_{ExtSGD}). \label{LBExtSGD}
\end{equation}
Notice that the construction of $\eta_t$ for algorithms in ${\cal A}_{ExtSGD}$ does \textit{not} depend on knowledge of the stochastic gradient $\nabla f(w_t;\xi_t)$. So, we do not consider step sizes that are adaptively computed based on $\nabla f(w_t;\xi_t)$.

As a disclaimer we notice that for some objective functions $F\in {\cal F}^N_{sm}[\mu,L]$ the expected convergence rate can be much better than what is stated in (\ref{LBExtSGD}); this is because $\gamma_t(\{F\}, {\cal A}_{ExtSGD})$ can be much smaller than  $\gamma_t({\cal F}^N_{sm}[\mu,L], {\cal A}_{ExtSGD})$, see (\ref{gammaineq}).
This is due to the specific nature of the objective function $F$ itself. However, without knowledge about this nature, one can only prove a general upper bound on the expected convergence rate $Y_t$ and any such upper bound must be at least the lower bound (\ref{LBExtSGD}).

Results (\ref{eqLB}) and (\ref{Yt}) of the previous section combined with (\ref{FAstoch}), (\ref{FAExtSGD}), and (\ref{LBExtSGD}) yield
\begin{eqnarray}
\frac{\log(2/\sqrt{e})}{432\cdot d}\frac{N}{\mu^2 t} \leq \gamma_t({\cal F}^N_{lip}[\mu,K,d], {\cal A}_{stoch}) &\leq  & \gamma_t({\cal F}^N_{str}[\mu], {\cal A}_{SGD}), \label{In1}\\
 \frac{1}{2}\frac{N}{\mu^2 t}(1- {\cal O}((\ln t)/t)) \leq \gamma_t({\cal F}^N_{esm}[\mu,L], {\cal A}_{ExtSGD}) &\leq &  \gamma_t({\cal F}^N_{str}[\mu], {\cal A}_{SGD}), \label{In2}\\
  \frac{1}{2}\frac{N}{\mu^2 t}(1- {\cal O}((\ln t)/t)) \leq
  \gamma_t({\cal F}^N_{sm}[\mu,L], {\cal A}_{ExtSGD}) 
  &\leq&  \gamma_t({\cal F}^N_{esm}[\mu,L], {\cal A}_{SGD}) \nonumber \\
  &\leq& \frac{16N}{\mu^2 t}(1+{\cal O}(1/t)). \label{In3}
\end{eqnarray}


We conclude the following observations (our contributions):
\begin{enumerate}
    \item 
The first inequality (\ref{In1}) is from \cite{AgarwalBartlettRavikumarEtAl}. Comparing  (\ref{In2}) to (\ref{In1}) shows that as a lower bound for $\gamma_t({\cal F}^N_{str}[\mu], {\cal A}_{SGD})$ (SGD for the class of strongly convex objective functions) our lower bound (\ref{LBExtSGD}) is dimension independent and improves the lower bound (\ref{eqLB}) of \cite{AgarwalBartlettRavikumarEtAl} by a factor $775\cdot d$. This is a significant improvement. 
\item However, our lower bound does not hold for the larger class ${\cal A}_{stoch}$. This teaches us that if we wish to reach smaller (better) expected convergence rates, then one approach is to step beyond SGD where our lower bound does not hold implying that within ${\cal A}_{stoch}$ there may be an opportunity to find an algorithm leading to at most a factor $32\cdot 775\cdot d$ smaller expected convergence rate compared to upper bound (\ref{In3}). This is the first exact quantification into what extent a better (practical) algorithm when compared to classical SGD can be found. 
E.g., Adam~\cite{kingma2014adam}, AdaGrad~\cite{duchi2011adaptive}, SGD-Momentum \cite{sutskever2013importance}, RMSProp \cite{zou2018sufficient} 
are all in ${\cal A}_{stoch}$ and can beat classical SGD by at most a 
factor $32\cdot 775\cdot d$. 
\item When searching for a better algorithm in ${\cal A}_{stoch}$ which significantly improves over SGD, it does not help to take an SGD-like algorithm which uses step sizes that are a function of iteratively computed estimates of $\nabla F(w_t)$ and $Y_t$ as this would keep such an algorithm in ${\cal A}_{ExtSGD}$ for which our lower bound is tight.
\item Another approach to reach smaller expected convergence rates is to stick with SGD but consider a smaller restricted class of objective functions for which  more/other information in the form of extra global parameters is available for adaptively computing $\eta_t$.
\item For strongly convex and smooth, respectively expected smooth, objective functions the algorithm $A\in {\cal A}_{SGD}$ with stepsizes $\eta_t = \frac{2}{\mu t+4L}$, respectively $\eta_t = \frac{2t+1}{(t+1)^2\mu}$ for $t> \frac{4L}{\mu}$ and $\eta_t=\frac{1}{2L}$ for $t\leq \frac{4L}{\mu}$, realizes the upper bound in (\ref{In3}) for all $t$. Inequalities (\ref{In3}) show that this algorithm is close to optimal: For each $t$, the best sequence of diminishing step sizes which minimizes $Y_t$ can at most achieve a constant (dimension independent) factor $32$ smaller expected convergence rate. 
\end{enumerate}

\section{Lower Bound for Extended SGD}
\label{sec:tight1}

In order to prove a lower bound
we propose a specific subclass of strongly convex and smooth objective functions $F$ and we show in the extended SGD setting how, based on recurrence (\ref{EqY}), to compute the \textit{optimal} step size $\eta_t$ as a function of $\mu$ and  $L$ and an oracle ${\cal U}$ with access to $N$, $\nabla F(w_t)$, and $Y_t$, i.e., this step size achieves the smallest
$Y_{t+1}$ at the $t$-th iteration.


We consider the following class of objective functions $F$:
We consider a multivariate normal distribution of a $d$-dimensional random vector $\xi$, i.e., $\xi \sim \mathcal{N}(m,\Sigma)$, where 
$m=\mathbb{E}[\xi]$ and $\Sigma = \mathbb{E}[(\xi-m)(\xi-m)^{\mathrm{T}}]$ 
is the (symmetric positive semi-definite) covariance matrix. The density function of $\xi$ is chosen as
$$
g(\xi) = \frac{\exp( \frac{- (\xi-m)^\mathrm{T}\Sigma^{-1}(\xi-m)}{2} )}{\sqrt{(2\pi)^d|\Sigma|}}. 
$$

We select  component functions $f(w;\xi)=s(\xi)\frac{\|w-\xi\|^2}{2}$,
where function $s(\xi)$ is constructed {\em a-priori} according to the following random process: 
 \begin{itemize} 
 \item With probability $1-\mu/L$, we draw $s(\xi)$ from the uniform distribution over interval $[0, \mu/(1-\mu/L)]$.
 \item With probability $\mu/L$, we draw $s(\xi)$ from the uniform distribution over interval $[0, L]$.
 \end{itemize}

The following theorem analyses the sequence of optimal step sizes for our class of objective functions and gives a lower bound on the corresponding expected convergence rates. The theorem states that we cannot find a better sequence of step sizes. In other words without any more additional information about the objective function (beyond  $\mu,L,N,Y_0,\dotsc,Y_t$ for computing $\eta_t$), we can at best prove a general upper bound which is at least the lower bound as stated in the theorem. 
The proof of the lower bound  is presented in the supplementary material:

\begin{thm}
\label{theorem:example_tightness}
We assume that component functions $f(w;\xi)$ are constructed according to the recipe described above with $\mu < L/18$. Then, the corresponding objective function is $\mu$-strongly convex and the component functions are $L$-smooth and convex. 
 
If we run Algorithm~\ref{sgd_algorithm} and assume that access to an oracle ${\cal U}$ with access to $N$, $\nabla F(w_t)$, and $Y_t$ is given at the $t$-th iteration (our extended SGD problem setting), then an exact expression for the optimal sequence of stepsizes $\eta_t$ 
based on $\mu,L,N, Y_0,\dotsc,Y_t$ can be given, i.e.,  this sequence of  stepsizes achieves the smallest possible $Y_{t+1}$ at the $t$-th iteration for all $t$. 
For this sequence of stepsizes,
\begin{equation}
 \label{tightLOW}
  Y_t \geq \frac{N}{2\mu} \frac{1}{\mu t+2\mu \ln(t+1) + W},
\end{equation}
where
\begin{align*}
    W= \frac{ L^2}{12(L-\mu)}. 
\end{align*}
\end{thm}

In the supplementary material we show numerical experiments in agreement with the presented theorem. 

\section{Conclusion}
\label{sec:conclusion}
We have studied the convergence of SGD by introducing a framework for comparing upper bounds and lower bounds and by proving a new lower bound based on straightforward calculus. The new lower bound is dimension independent and improves a factor $775\cdot d$ over previous work \cite{AgarwalBartlettRavikumarEtAl} applied to SGD, shows the optimality of step sizes in \cite{Nguyen2018,gower2019sgd},
and shows that practical algorithms like Adam~\cite{kingma2014adam}, AdaGrad~\cite{duchi2011adaptive}, SGD-Momentum \cite{sutskever2013importance}, RMSProp \cite{zou2018sufficient} for strongly convex objective functions
can at most achieve a factor $32\cdot 775\cdot d$ smaller expected convergence rate compared to classical SGD.

\section*{Acknowledgement}
We thank the reviewers for useful suggestions to improve the paper. Phuong Ha Nguyen and Marten van Dijk were supported in part by AFOSR MURI under award number FA9550-14-1-0351. 


\bibliography{reference}
\bibliographystyle{plainnat}


\appendix
\newpage

\section*{Supplementary Material}


%


\section{Proof}
\label{supproof}

We extend Theorem~\ref{theorem:example_tightness} with an upper bound used in our numerical experiments.

\textbf{Theorem~\ref{theorem:example_tightness}}
{\em We assume that component functions $f(w;\xi)$ are constructed according to the recipe described in Section \ref{sec:tight1} with $\mu < L/18$. Then, the corresponding objective function is $\mu$-strongly convex and the component functions are $L$-smooth and convex. 

If we run Algorithm~\ref{sgd_algorithm} and assume that access to an oracle ${\cal U}$ with access to $N$, $\nabla F(w_t)$, and $Y_t$ is given at the $t$-th iteration (our extended SGD problem setting), then an exact expression for the optimal sequence of stepsizes $\eta_t$ 
based on $\mu,L,N, Y_0,\dotsc,Y_t$ can be given, i.e.,  this sequence of  stepsizes achieves the smallest possible $Y_{t+1}$ at the $t$-th iteration for all $t$. 
For this
sequence of stepsizes,}
\begin{equation*}
  Y_t \geq \frac{N}{2\mu} \frac{1}{\mu t+2\mu \ln(t+1) + W},
\end{equation*}
where $W= \frac{ L^2}{12(L-\mu)}$
and for $t \geq T' = \frac{20 L}{\mu}$,
\begin{equation}
 \label{tightUP}
  Y_t \leq \frac{16N}{\mu}\frac{1}{\mu t-16 L}.
\end{equation}

\begin{proof}
We first restrict oracle ${\cal U}$ to only supply information about $N$ and $Y_t$ at the $t$-th iteration. At the end of this proof we show that our arguments generalize to the more powerful oracle ${\cal U}$ which also provides the full gradient $\nabla F(w_t)$  at the $t$-th iteration.

Clearly, $f(w;\xi)$ is $s(\xi)$-smooth where the maximum value of $s(\xi)$ is equal to $L$. That is, all functions $f(w;\xi)$ are $L$-smooth (and we cannot claim a smaller smoothness parameter). We notice that

 $$ \mathbb{E}_{\xi}[s(\xi)]  = (1-\mu/L) \frac{\mu/(1-\mu/L)}{2} + (\mu/L) \frac{L}{2} = \mu$$
 and
 \begin{align*}
  \mathbb{E}_{\xi}[s(\xi)^2]  &=  (1-\mu/L) \frac{(\mu/(1-\mu/L))^2}{12}+ (\mu/L) \frac{L^2}{12}\\&= \frac{\mu (L + \frac{\mu}{1-\mu/L})}{12}= \frac{\mu L^2}{12(L-\mu)}.  
 \end{align*}
 
 
With respect to $f(w;\xi)$ and distribution $g(\xi)$ we define
  $$F(w) = \mathbb{E}_{\xi}[f(w;\xi)] = \mathbb{E}_{\xi}[s(\xi) \frac{\|w-\xi\|^2}{2}].$$
  \textit{Since $s(\xi)$ only assigns a random variable to $\xi$ which is drawn from a distribution whose description is not a function of $\xi$, random variables $s(\xi)$ and $\xi$ are statistically independent}. Therefore, $F(w)=$
$$ \mathbb{E}_{\xi}[s(\xi) \frac{\|w-\xi\|^2}{2}]= \mathbb{E}_{\xi}[s(\xi)]  \mathbb{E}_{\xi}[\frac{\|w-\xi\|^2}{2}] = \mu  \mathbb{E}_{\xi}[\frac{\|w-\xi\|^2}{2}] $$
  
 Notice:
\begin{enumerate}
\item $\|w-\xi\|^2 = \|(w-m)+(m-\xi)\|^2= \|w-m\|^2 + 2 \langle w-m,m-\xi \rangle + \|m-\xi \|^2$. 
\item Since $m=\mathbb{E}[ \xi ]$, we have $\mathbb{E}[ m-\xi ] =0$.
\item   $\mathbb{E}[ \|m-\xi\|^2 ]=\sum_{i=1}^d \mathbb{E}[ (m_i-\xi_i)^2 ] = \sum_{i=1}^d \Sigma_{i,i} = \mathrm{Tr}(\Sigma)$. 
\end{enumerate}
 Therefore, 
 $F(w) = \mu  \mathbb{E}_{\xi}[\frac{\|w-\xi\|^2}{2}] = \mu \frac{\|w-m\|^2}{2} + \mu \frac{\mathrm{Tr}(\Sigma)}{2},$ and this shows $F$ is $\mu$-strongly convex and has minimum $w_* = m$.
 
 Since
\begin{align*}
 \nabla_w [ \|w-\xi\|^2 ] &= \nabla_w [ \langle w,w \rangle - 2\langle w,\xi \rangle + \langle \xi,\xi \rangle ]\\& = 2w -2 \xi = 2(w-\xi),   
\end{align*} 
 

we have $$ \nabla_w f(w;\xi)= s(\xi)(w-\xi). $$
In our notation 
$$N=2 \mathbb{E}_\xi[\|\nabla f(w_*;\xi) \|^2]=2 \mathbb{E}_{\xi}[s(\xi)^2 \|w_* -\xi \|^2].$$
By using similar arguments as used above we can split the expectation and obtain
$$ N = 2 \mathbb{E}_{\xi}[s(\xi)^2 \|w_* -\xi \|^2] = 2  \mathbb{E}_{\xi}[s(\xi)^2] \mathbb{E}_{\xi}[ \|w_* -\xi \|^2].$$
We already calculated ($w_*=m$)
$$\mathbb{E}_{\xi}[ \|w_* -\xi \|^2] = \|w_*-m\|^2 + \mathrm{Tr}(\Sigma)= \mathrm{Tr}(\Sigma)$$
and we know $$ \mathbb{E}_{\xi}[s(\xi)^2] = \frac{\mu L^2}{12(L-\mu)}.$$ 

This yields
$$ N =  2  \mathbb{E}_{\xi}[s(\xi)^2] \mathbb{E}_{\xi}[ \|w_* -\xi \|^2]
= \frac{\mu L^2}{6(L-\mu)} \mathrm{Tr}(\Sigma).$$

In the SGD algorithm we compute
\begin{align*}
w_{t+1} &= w_t - \eta_t \nabla f(w_t;\xi_t) \\ &= w_t - \eta_t s(\xi_t) (w_t - \xi_t) \\&= (1-\eta_t s(\xi_t)) w_t + \eta_t s(\xi_t) \xi_t. 
\end{align*}


We draw $\xi$ from its distribution and set $w_0 = \xi$. Therefore,
$$ Y_0 = \mathbb{E}[ \|w_0-w_*\|^2] = \mathbb{E}[ \|\xi-w_*\|^2 ] = \mathrm{Tr}(\Sigma).$$


Let $\mathcal{F}_t=\sigma(w_0,\xi_0,\dotsc,\xi_{t-1})$ be the $\sigma$-algebra generated by $w_0,\xi_0,\dotsc,\xi_{t-1}$. We derive $\mathbb{E}[ \|w_{t+1}-w_*\|^2 | \mathcal{F}_t ]$
$$
=\mathbb{E}[ \|(1-\eta_t s(\xi_t)) (w_t -w_*) + \eta_t s(\xi_t) (\xi_t - w_*) \|^2 | \mathcal{F}_t ]
$$
which is equal to 
\begin{align}
&\mathbb{E}[ (1-\eta_t s(\xi_t))^2 \|w_t -w_*\|^2 \nonumber \\ 
& + 2\eta_t s(\xi_t)(1-\eta_t s(\xi_t) ) \langle  w_t-w_*, \xi_t - w_* \rangle \nonumber \\ 
&+  \eta_t^2 s(\xi_t)^2 \|\xi_t -w_*\|^2 | \mathcal{F}_t].  \label{eq:xxx1}   
\end{align}


Given $\mathcal{F}_t$, $w_t$ is not a random variable. Furthermore, we can use linearity of taking expectations and as above split expectations:
\begin{align}
&\mathbb{E}[(1-\eta_t s(\xi_t))^2] \|w_t -w_*\|^2 \nonumber \\
&+ \mathbb{E}[2\eta_t s(\xi_t)(1-\eta_t s(\xi_t) )] \langle w_t-w_*, \mathbb{E}[\xi_t - w_*] \rangle \nonumber \\
&+ \mathbb{E}[ \eta_t^2 s(\xi_t)^2] \mathbb{E}[ \|\xi_t -w_*\|^2]. \label{eq:xxx2}    
\end{align}


Again notice that $\mathbb{E}[\xi_t -w_*] = 0$ and $\mathbb{E}[ \| \xi_t - w_* \|^2 ] = \mathrm{Tr}(\Sigma)$. So, 
 $\mathbb{E}[\|w_{t+1}-w_*\|^2 | \mathcal{F}_t ]$ 
is equal to
\begin{eqnarray*}
 &&
\mathbb{E}[(1-\eta_t s(\xi_t))^2] \|w_t -w_*\|^2  +\eta_t^2 \frac{N}{2}  \\
&=& (1- 2\eta_t \mu + \eta_t^2 \frac{\mu L^2 }{12(L-\mu)}) \|w_t -w_*\|^2  +\eta_t^2 \frac{N}{2} \\
&=& (1-\mu \eta_t (2 - \frac{\eta_t}{12} \frac{L^2}{L-\mu}))  \|w_t -w_*\|^2  +\eta_t^2 \frac{N}{2} .
\end{eqnarray*}
In terms of $Y_t=\mathbb{E}[\|w_{t}-w_*\|^2]$, by taking the full expectation (also over $\mathcal{F}_t$) we get
\begin{equation} Y_{t+1} = (1-\mu \eta_t (2 - \frac{\eta_t}{12} \frac{L^2}{L-\mu}))  Y_t  +\eta_t^2 \frac{N}{2}. \label{eq:Yeta} \end{equation}
This is very close to recurrence (\ref{main_ineq_sgd_new02}).

Equation (\ref{eq:Yeta}) expresses $Y_{t+1}$ as a function $Y_{t+1}(\eta_t, Y_t)$ of $\eta_t$ and $Y_t$. Given $Y_0$, we want to minimize $Y_{t+1}$ with respect to the step sizes $\eta_t, \eta_{t-1}, \ldots, \eta_0$. For $i<t$ we derive
$$ \frac{\partial  Y_{t+1}}{\partial \eta_i}
= \frac{\partial  Y_{t+1}}{\partial Y_t} \frac{\partial Y_t}{\partial \eta_i} =
(1-\mu \eta_t (2 - \frac{\eta_t}{12} \frac{L^2}{L-\mu})) 
\frac{\partial Y_t}{\partial \eta_i} $$
and
for $i=t$ we derive 
\begin{equation}
 \label{optimalstepsize}
 \frac{\partial  Y_{t+1}}{\partial \eta_i}
= -2\mu Y_t + 2\mu  \frac{\eta_t}{12} \frac{L^2}{L-\mu}  Y_t + N\eta_t.
 \end{equation}
 We reach a stationary point for $Y_{t+1}$ as a function of step sizes $\eta_t, \eta_{t-1}, \ldots, \eta_0$ if each of the partial derivatives with respect to $\eta_i$ is equal to $0$. We notice that if for all $t$
 \begin{equation} 1-\mu \eta_t (2 - \frac{\eta_t}{12} \frac{L^2}{L-\mu}) >0,\label{Eqzero} \end{equation}
 then, for $i<t$, $\frac{\partial  Y_{t+1}}{\partial \eta_i}=0$ if and only if $\frac{\partial Y_t}{\partial \eta_i} =0$. This implies that $Y_{t+1}$ has a stationary point if and only if 
 $$\forall_{0\leq i\leq t} \ 
 \frac{\partial  Y_{i+1}}{\partial \eta_i}=0.$$
 Hence, if a step size sequence satisfies this for all $t$, then it leads to stationary points for all $Y_{t+1}$ as function of $\eta_t, \eta_{t-1}, \ldots, \eta_0$. So, such a sequence of step sizes simultaneously achieves stationary points for all $Y_{t+1}$. 
 
For the argument to hold, we need to prove (\ref{Eqzero}). 
The left hand side of (\ref{Eqzero}) achieves its minimum value $$1-12\mu\frac{L-\mu}{L^2}$$ for $\eta_t=12\frac{L-\mu}{L^2}$. 
For $\mu < \frac{L}{12}$, $12\mu(L-\mu)<12\mu L < L^2$ implying that this minimum value is larger than zero.


 As explained above the optimal step size $\eta_t$ in a sequence of optimal step sizes that minimizes all expected convergence rates $Y_t$ is computed by taking the derivative of $Y_{t+1}$ with respect to $\eta_t$. This derivative is equal to (\ref{optimalstepsize}) and shows that
  the minimum is achieved for
 \begin{equation}
 \label{optimal_learningrate}
 \eta_t = \frac{ 2\mu Y_t }{N + \frac{\mu L^2}{6(L-\mu)}  Y_t}
 \end{equation}
 giving, see (\ref{eq:Yeta}),
 \begin{align}
Y_{t+1} &=Y_t - \frac{ 2\mu^2 Y_t^2 }{N + \frac{\mu L^2}{6(L-\mu)}  Y_t} \nonumber \\
  &= Y_t - \frac{ 2\mu^2 Y_t^2 }{N (1+ Y_t/ \mathrm{Tr}(\Sigma))}. \label{ytoptimal}
 \end{align}

  We note that $Y_{t+1} \leq Y_t$ for any $t \geq 0$. We proceed by proving a lower bound on $Y_t$.
  Clearly,
  \begin{equation}
  \label{q1111}
  Y_{t+1} \geq Y_t - \frac{ 2\mu^2 Y_t^2 }{N}
  \end{equation}
  Let us define $\gamma = 2\mu^2/N$. We can rewrite (\ref{q1111}) as follows:
  \begin{eqnarray}
    \gamma Y_{t+1} &\geq& \gamma Y_t(1-\gamma Y_t) \text{ or} \nonumber \\
    (\gamma Y_{t+1})^{-1} &\leq& 1 + (\gamma Y_{t})^{-1} + \frac{1}{(\gamma Y_{t})^{-1}-1}. \label{eqYYYY}
  \end{eqnarray}

In order to make the inequality above correct, we require $1-\gamma Y_t > 0$ for any $t\geq 0$. Since $Y_{t+1} \leq Y_t$, we only need $Y_0 < \frac{1}{\gamma}$. This is implied by $Y_0=\mathrm{Tr}(\Sigma) < \frac{2}{3\gamma}$, a condition which is needed in the next sequence of arguments. This stronger condition means that we need $$\mathrm{Tr}(\Sigma) < \frac{N}{3\mu^2}, \ \text{i.e.,} \ \mathrm{Tr}(\Sigma)<\frac{\mu L^2}{6(L-\mu)} \frac{\mathrm{Tr}(\Sigma)}{3\mu^2}$$ after substituting $N$. This is equivalent to $\mu < \frac{L^2}{18(L-\mu)}$ which is true for $\mu < L/18$.    

By using induction on $t$, upper bound (\ref{eqYYYY}) implies 
\begin{equation}
(\gamma Y_{t+1})^{-1} \leq (t+1) + (\gamma Y_{0})^{-1} + \sum_{i=0}^t \frac{1}{(\gamma Y_{i})^{-1}-1}.    \label{eq:YYY}
\end{equation}
In order to further upper bound the sum in the right hand side, we first find a lower bound on $(\gamma Y_{i})^{-1}$. We rewrite equation (\ref{ytoptimal}) as
$$ (\gamma Y_{t+1}) = (\gamma Y_t) (1-\frac{(\gamma Y_t)}{1+ Y_t/\mathrm{Tr}(\Sigma)}).$$
Since $Y_t\leq Y_0=\mathrm{Tr}(\Sigma)$, we have
$$ (\gamma Y_{t+1}) \leq (\gamma Y_t) (1-\frac{(\gamma Y_t)}{2}).$$
This translates into
\begin{eqnarray*} (\gamma Y_{t+1})^{-1} &\geq& \frac{(\gamma Y_t)^{-1} }{1-(\gamma Y_t)/2} = 
\frac{(\gamma Y_t)^{-2} }{(\gamma Y_t)^{-1}-1/2}\\
&=& \frac{1}{2} + (\gamma Y_t)^{-1} + \frac{1}{4(\gamma Y_t)^{-1}-2} \\
&\geq & \frac{1}{2} + (\gamma Y_t)^{-1},
\end{eqnarray*}
where the last inequality follows from $(\gamma Y_t)^{-1}>(\gamma Y_0)^{-1}= (\gamma \mathrm{Tr}(\Sigma))^{-1} >1$ making $4(\gamma Y_t)^{-1}-2$ positive.

The resulting inequality leads to a recurrence and by using induction on $t$ we obtain
$$(\gamma Y_{t+1})^{-1}\geq (t+1)/2 + (\gamma Y_{0})^{-1}.$$
Now we are able to upper bound
\begin{align*}
\sum_{i=0}^t \frac{1}{(\gamma Y_{i})^{-1}-1} &\leq \sum_{i=0}^t \frac{1}{i/2+ (\gamma Y_{0})^{-1}-1}\\
&=2\sum_{i=0}^t \frac{1}{i+ 2((\gamma Y_{0})^{-1}-1)}. 
\end{align*}

We showed earlier that $\mu<L/18$ implies $Y_0<\frac{2}{3\gamma}$. Substituting this upper bound in our derivation leads to
$$\sum_{i=0}^t \frac{1}{(\gamma Y_{i})^{-1}-1} \leq 2 \sum_{i=0}^t \frac{1}{i+ 1} \leq 2\ln(t+2).$$
Combining with (\ref{eq:YYY}) we have the following inequality: 
\begin{align*}
(\gamma Y_{t+1})^{-1} &\leq (t+1) + (\gamma Y_0)^{-1} + 2\ln(t+2).
\end{align*}

Reordering, substituting $\gamma = 2\mu^2/N$, and replacing  $t+1$ by $t$ yields, for $t\geq 0$, the lower bound
\begin{eqnarray*}
 Y_t &\geq&  \frac{N}{2\mu} \frac{1}{\mu t + N/(2\mu  Y_0) + 2\mu\ln(t+1)} \nonumber \\
&=& \frac{N}{2\mu} \frac{1}{\mu t+2\mu \ln(t+1) + W},  
\end{eqnarray*}
where 
\begin{align*}
W=& N/(2\mu Y_0) =  \frac{L^2}{12(L-\mu)}.
\end{align*}

We now extend oracle ${\cal U}$ to also provide information about 
\textit{full gradient $\nabla F(w_t)$}  at the $t$-th iteration.
The above proof generalizes to this more powerful oracle. This is because of the reason why we are allowed to transform (\ref{eq:xxx1}) into (\ref{eq:xxx2}), i.e., $\eta_t$ and $\xi_t$ must be independent to get (\ref{eq:xxx2}) from (\ref{eq:xxx1}). If the construction of $\eta_t$ does not depend on $\xi_t$ (or $\nabla f(w_t;\xi_t)$), then only $Y_t$ is required to construct the optimal stepsize $\eta_t$. It implies that the information of $\nabla F(w_t)$ is not useful and we can borrow the above proof to arrive at the lower bound of this theorem.    

The upper bound for $Y_t$ comes from the following fact. If we run Algorithm~\ref{sgd_algorithm} with step size $\eta'_t=\frac{2}{\mu t +4L}$ for $t \geq 0$ in~\cite{Nguyen2018}, then we have from \cite{Nguyen2018}
an expected convergence rate
$$Y'_t \leq \frac{16N}{\mu}\frac{1}{\mu(t-T')+4L}$$ 
for $t\geq T'$, where $$T'=\frac{4L}{\mu}\max\{\frac{L\mu Y_0}{N},1\}-\frac{4L}{\mu}.$$ Substituting $$N = \frac{\mu L^2}{6(L-\mu)} \mathrm{Tr}(\Sigma) \ \text{and} \ Y_0=\mathrm{Tr}(\Sigma)$$ yields $T'\leq \frac{20L}{\mu}$. Since $\eta_t$ is the most optimal step size and $\eta'
_t$ is not, $Y_t \leq Y'_t$. I.e., we have for $t\geq \frac{20L}{\mu} \geq T'$,
$$Y_t \leq \frac{16N}{\mu}\frac{1}{\mu(t-\frac{20L}{\mu})+4L}
=
\frac{16N}{\mu}\frac{1}{\mu t-16 L}.$$
\end{proof}

%

\section{Numerical Experiments}
\label{sec_numerical_tightness}

We verify our theory by considering simulations with different values of sample size $n$ (1000, 10000, and 100000) and vector size $d$ (10, 100, and 1000). 
We generate $m \in \mathbb{R}^d$ and a diagonal matrix $\Sigma \in \mathbb{R}^{d \times d}$ by drawing each element in $m$ and each element on the diagonal of $\Sigma$ at random from a uniform distribution over $[0,1]$. We have $L=1$ and $\mu=1/n$ where $n$ is the number of samples. Hence the condition number $L/\mu$ is equal to $n$ and represents the number of SGD iterations in a single epoch. We experimented with 10 runs and reported the average results. 

\begin{figure*}[h]
 \centering
 \includegraphics[width=0.32\textwidth]{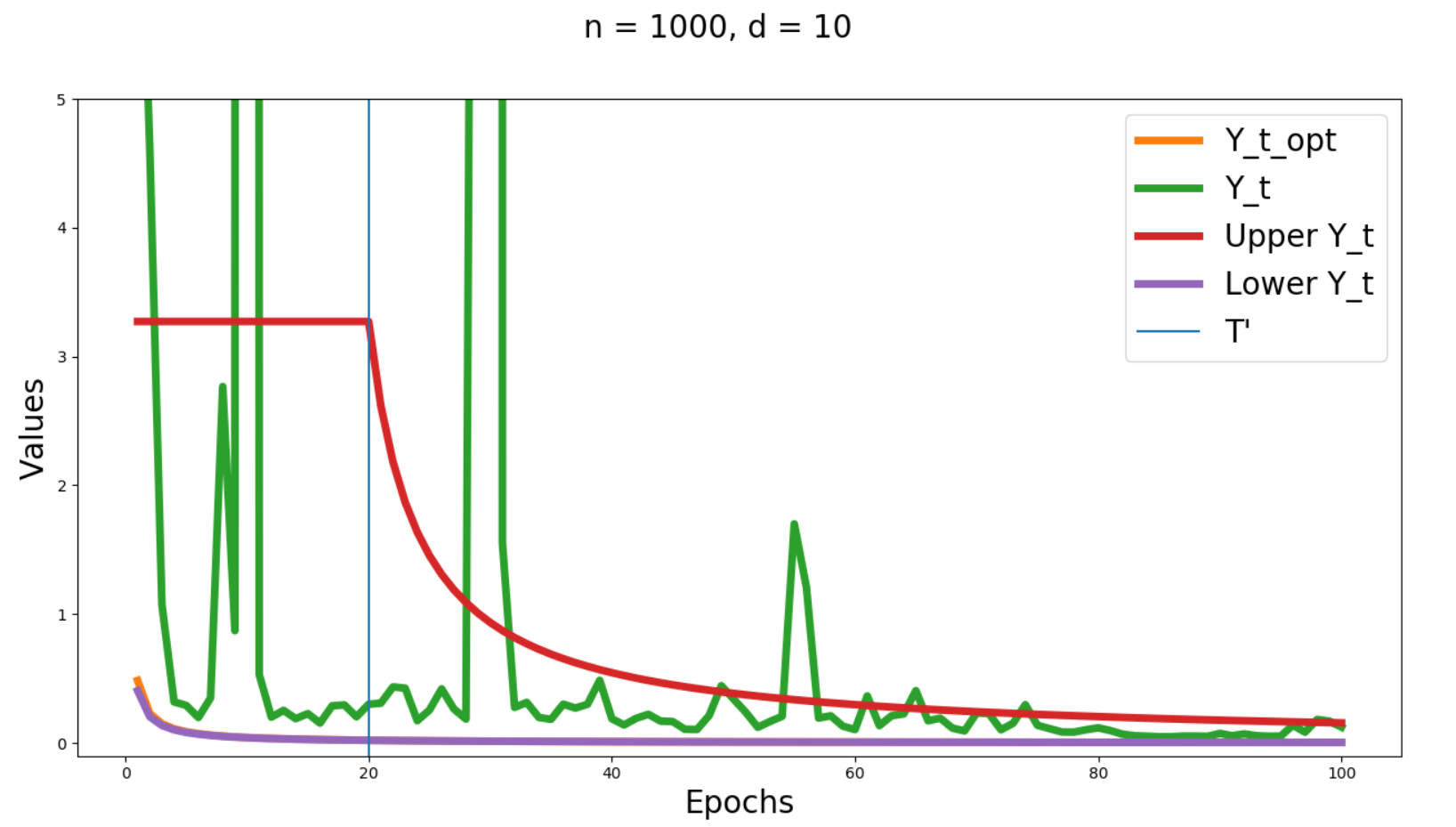}
 \includegraphics[width=0.32\textwidth]{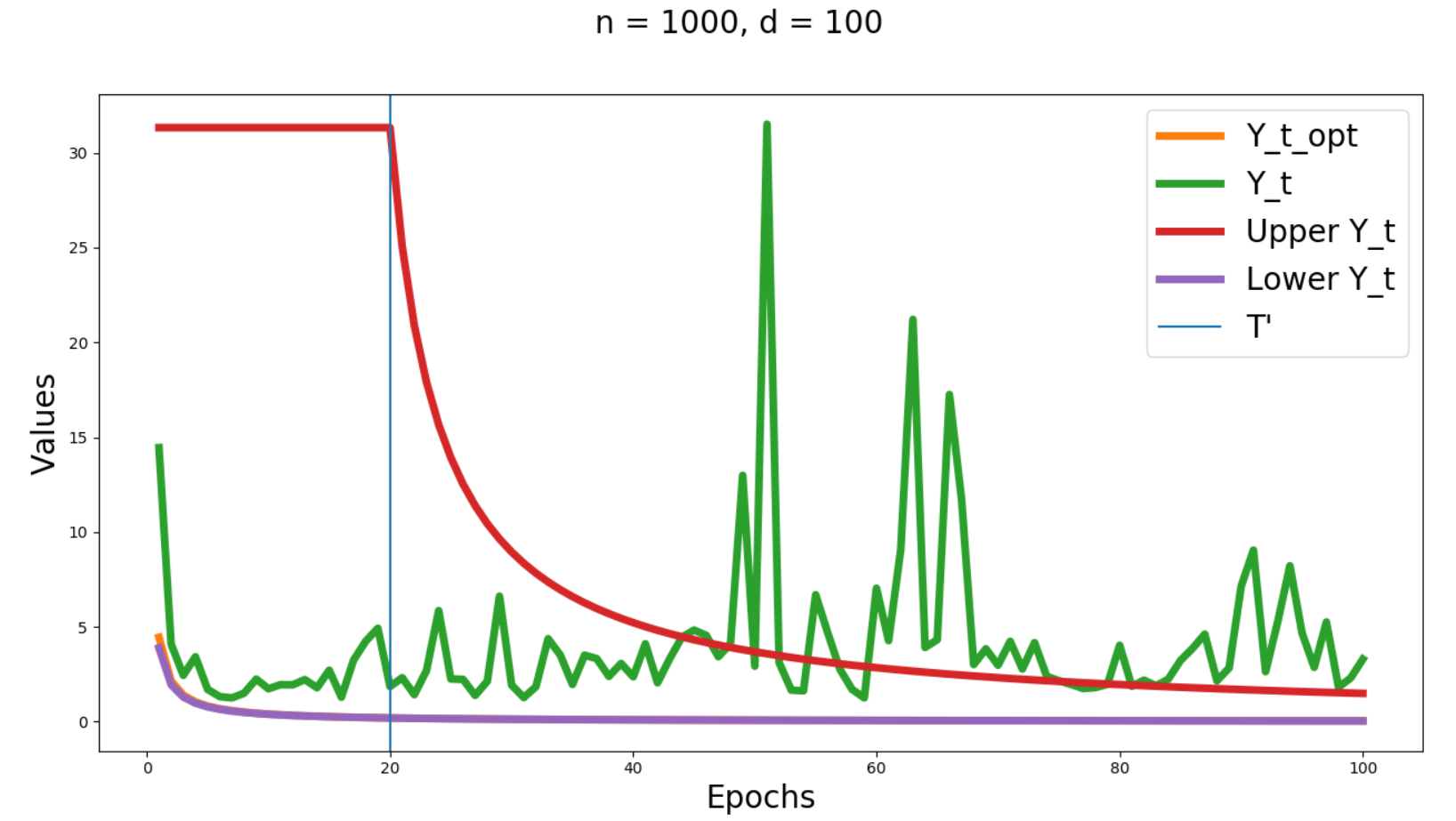}
 \includegraphics[width=0.32\textwidth]{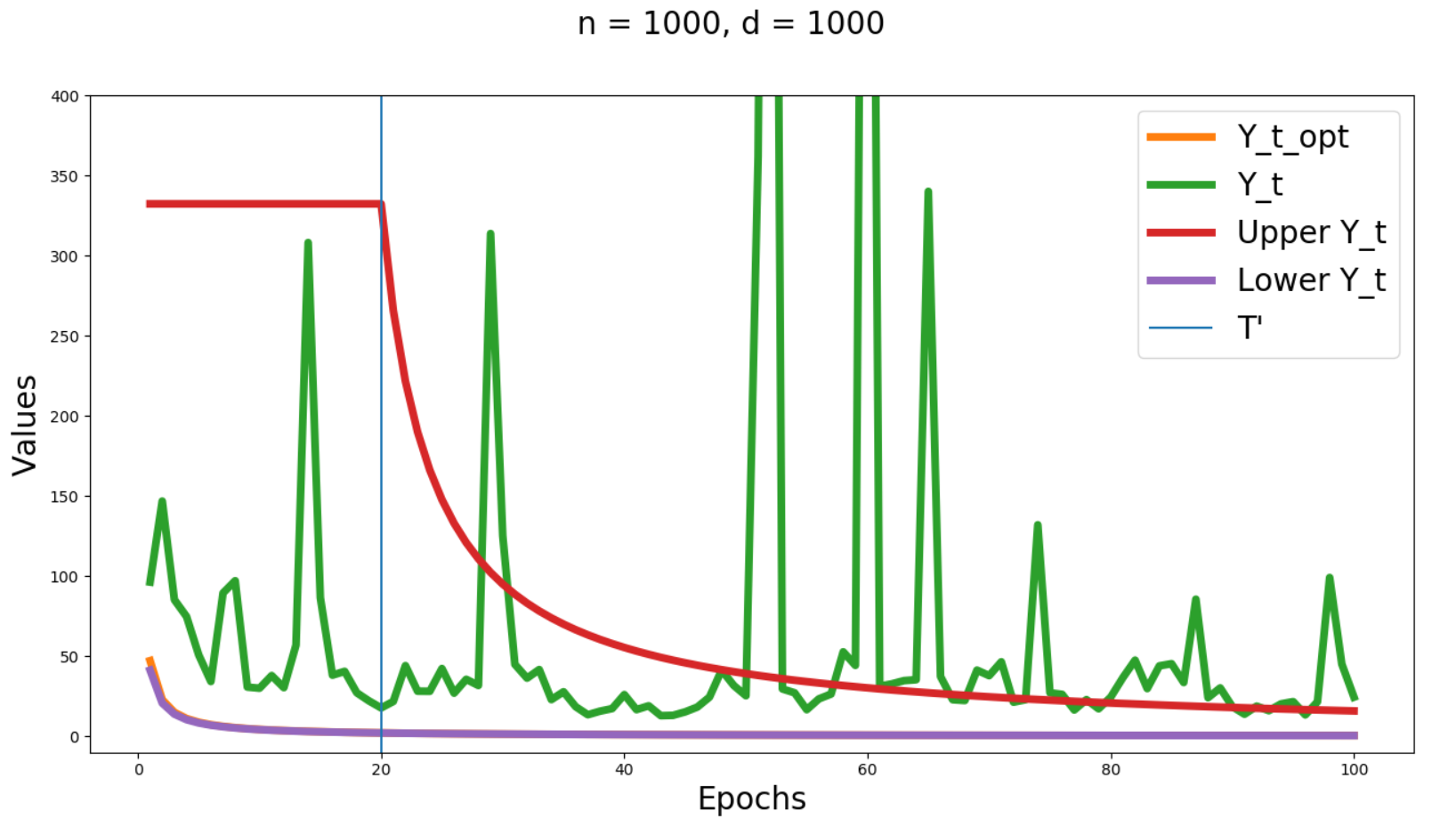}
 \includegraphics[width=0.32\textwidth]{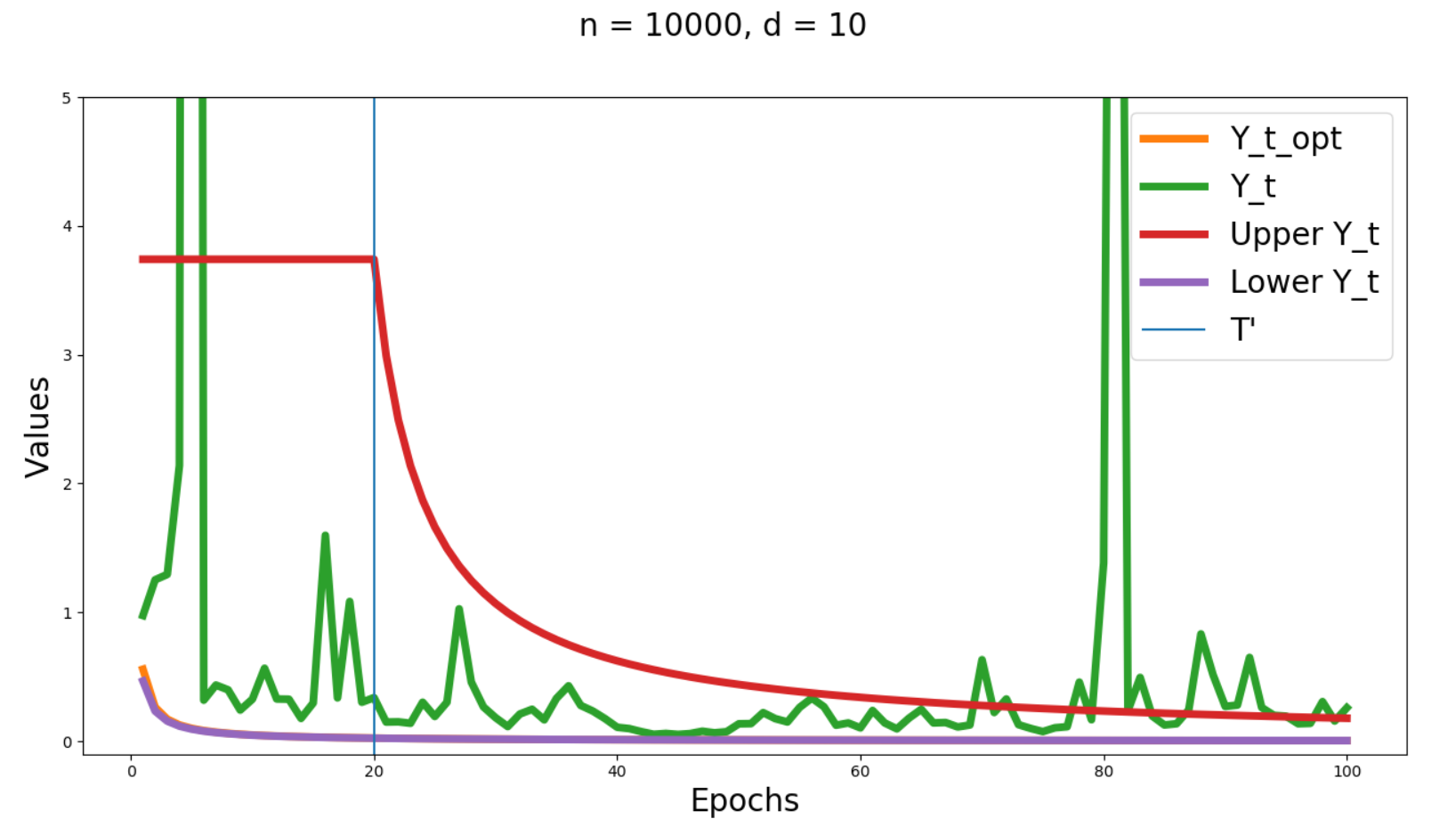}
 \includegraphics[width=0.32\textwidth]{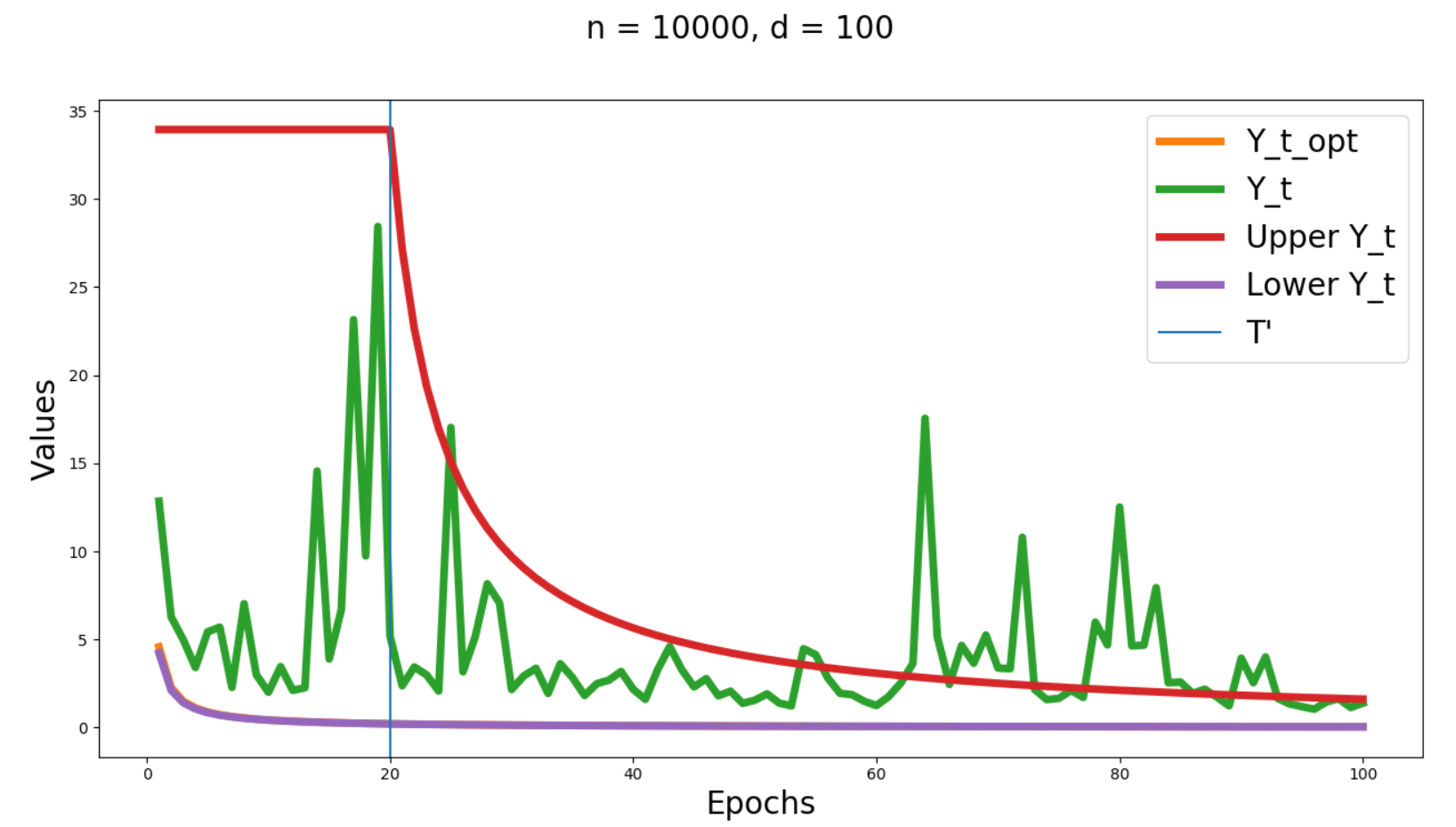}
 \includegraphics[width=0.32\textwidth]{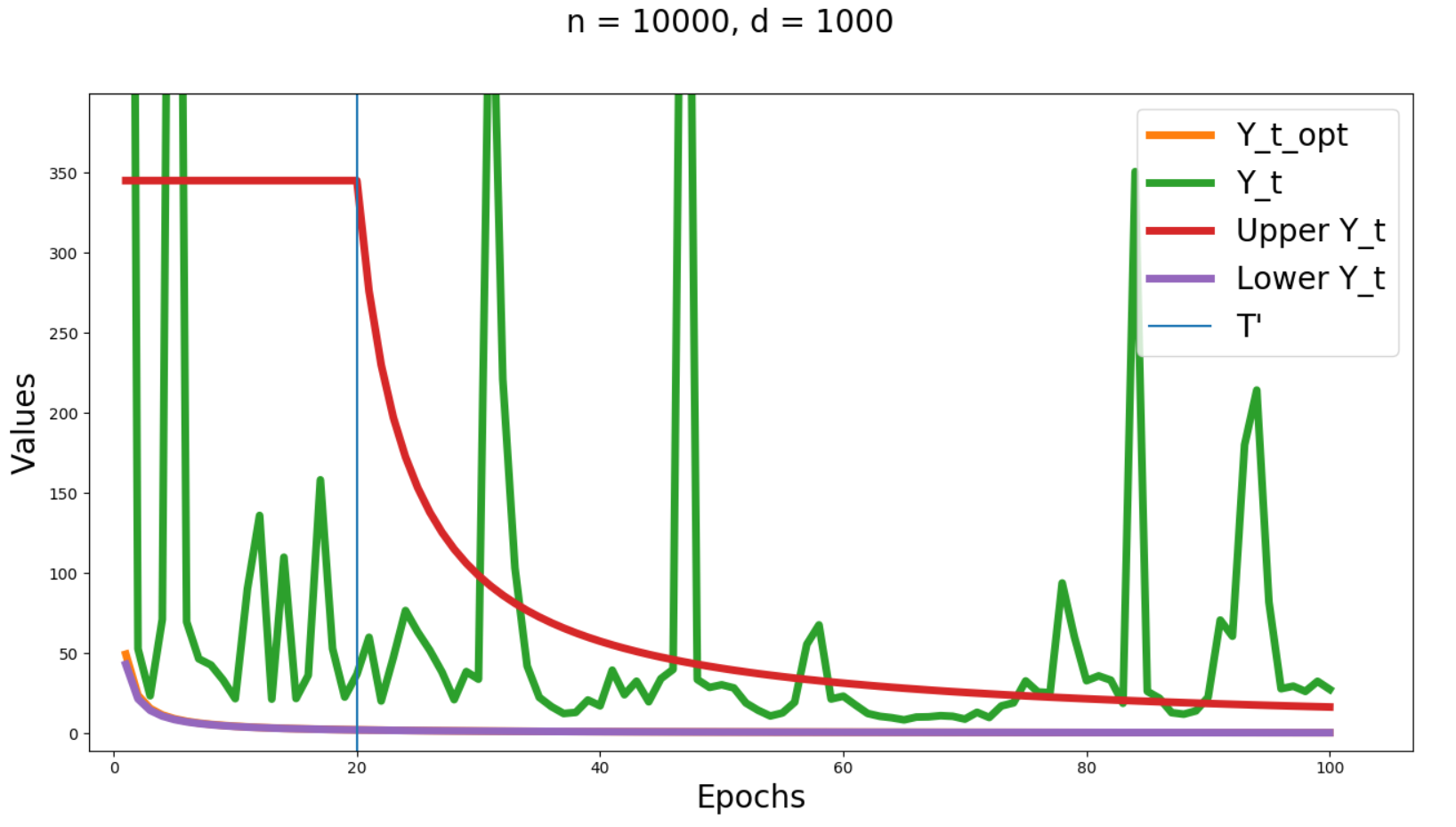}
 \includegraphics[width=0.32\textwidth]{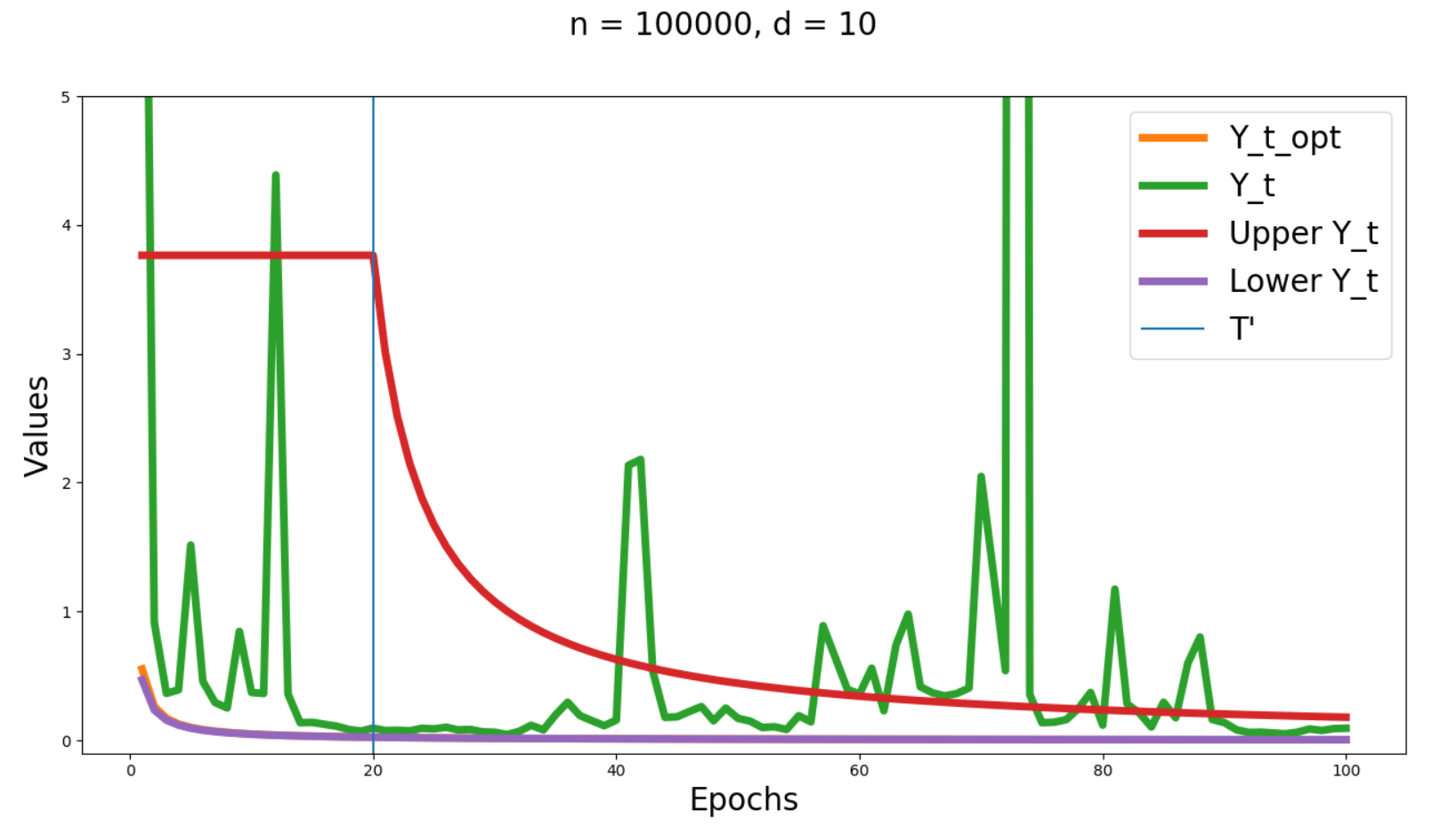}
 \includegraphics[width=0.32\textwidth]{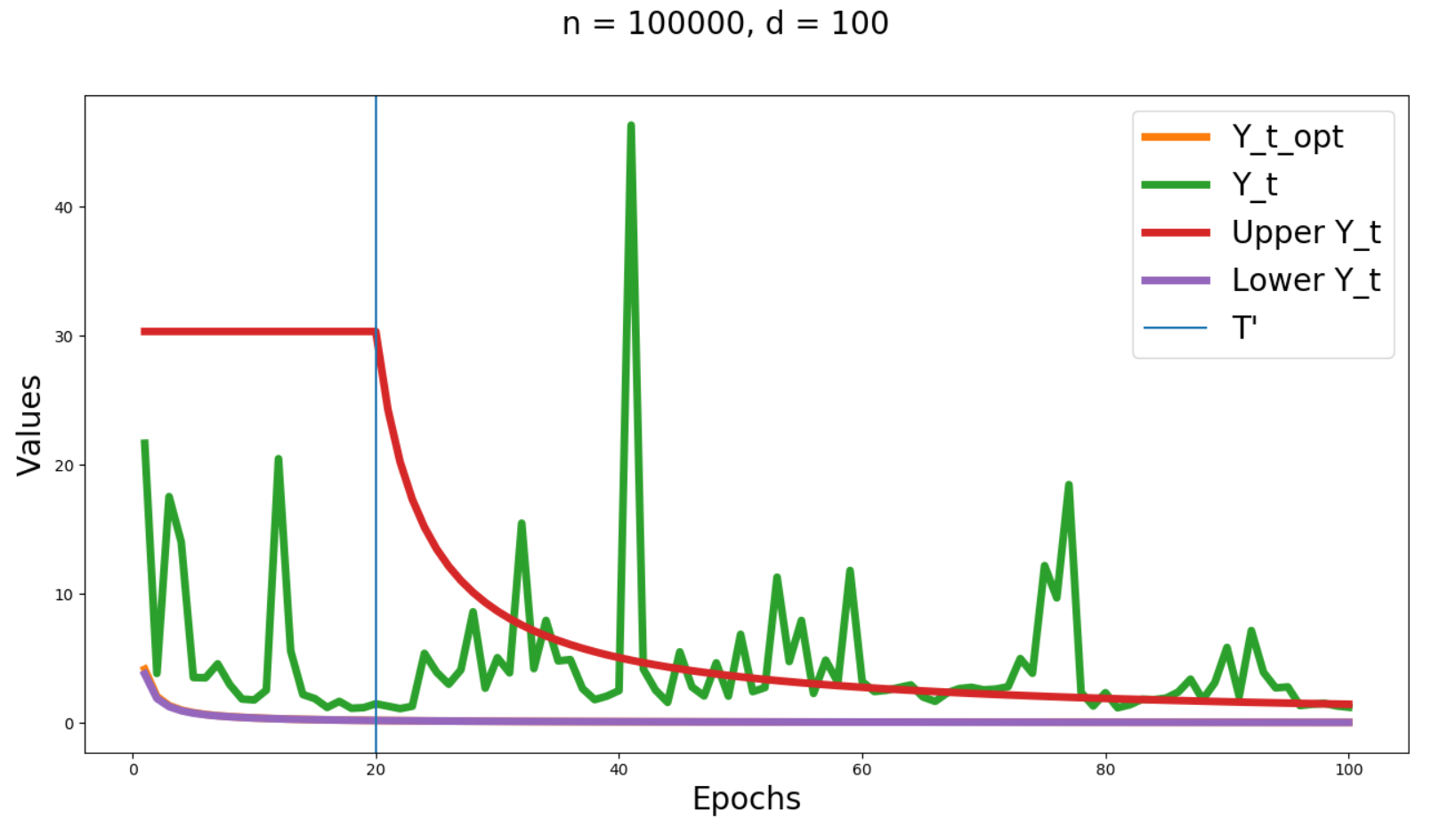}
 \includegraphics[width=0.32\textwidth]{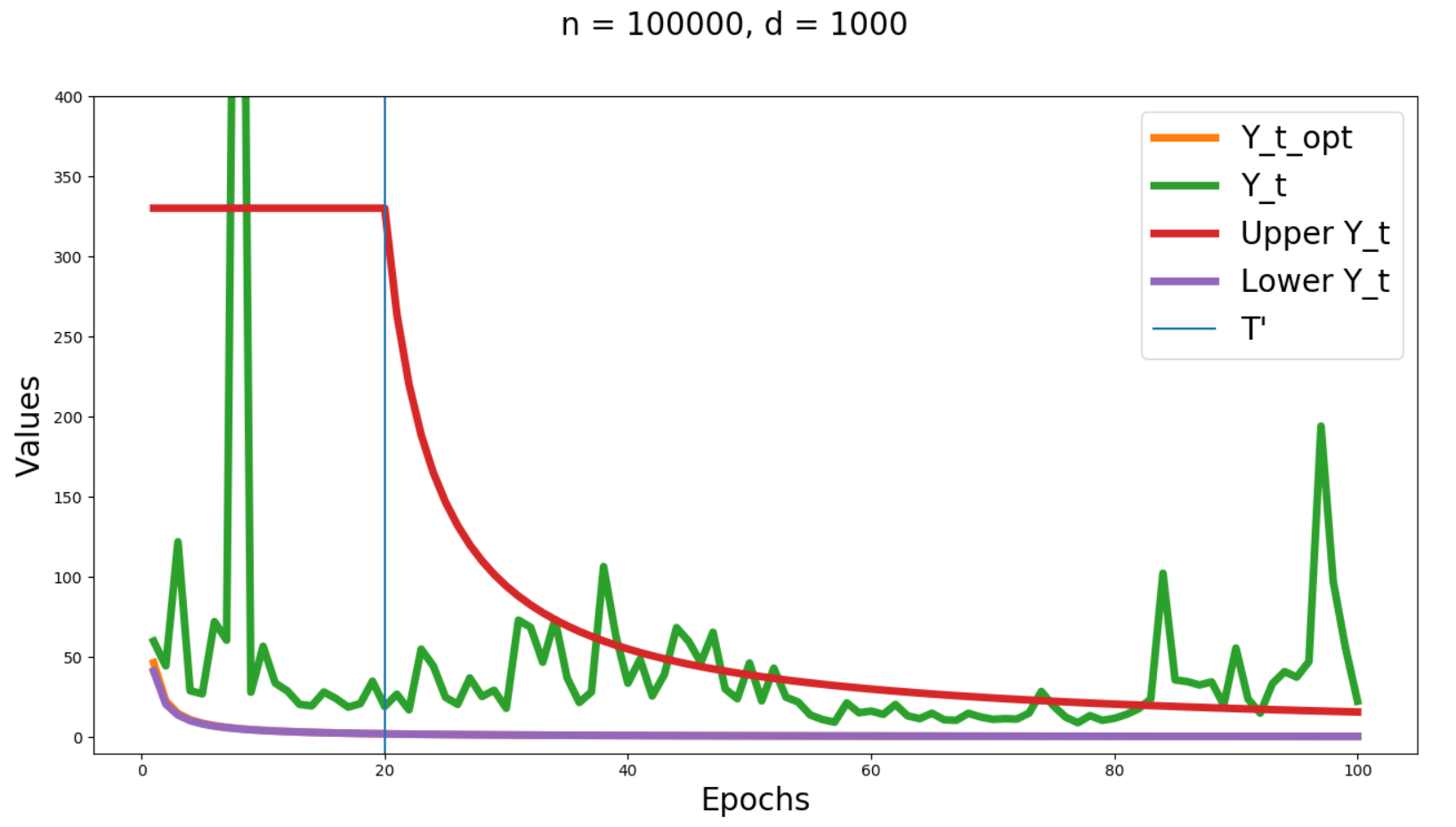} 
   \caption{$Y_t$ and its upper and lower bounds}
  \label{figure_n1}
 \end{figure*}

\begin{figure*}[!]
 \centering
 \includegraphics[width=0.32\textwidth]{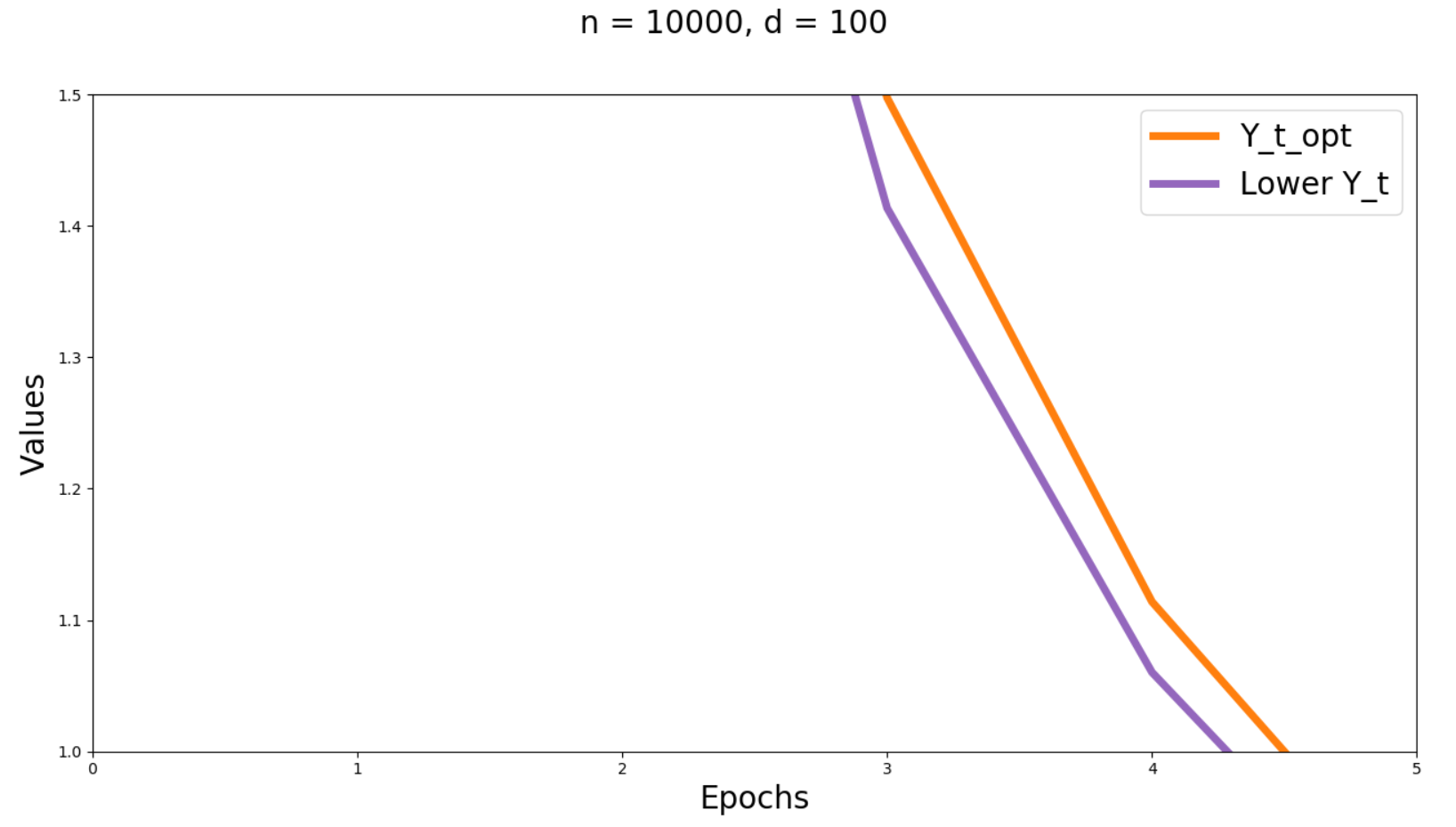}
 \includegraphics[width=0.32\textwidth]{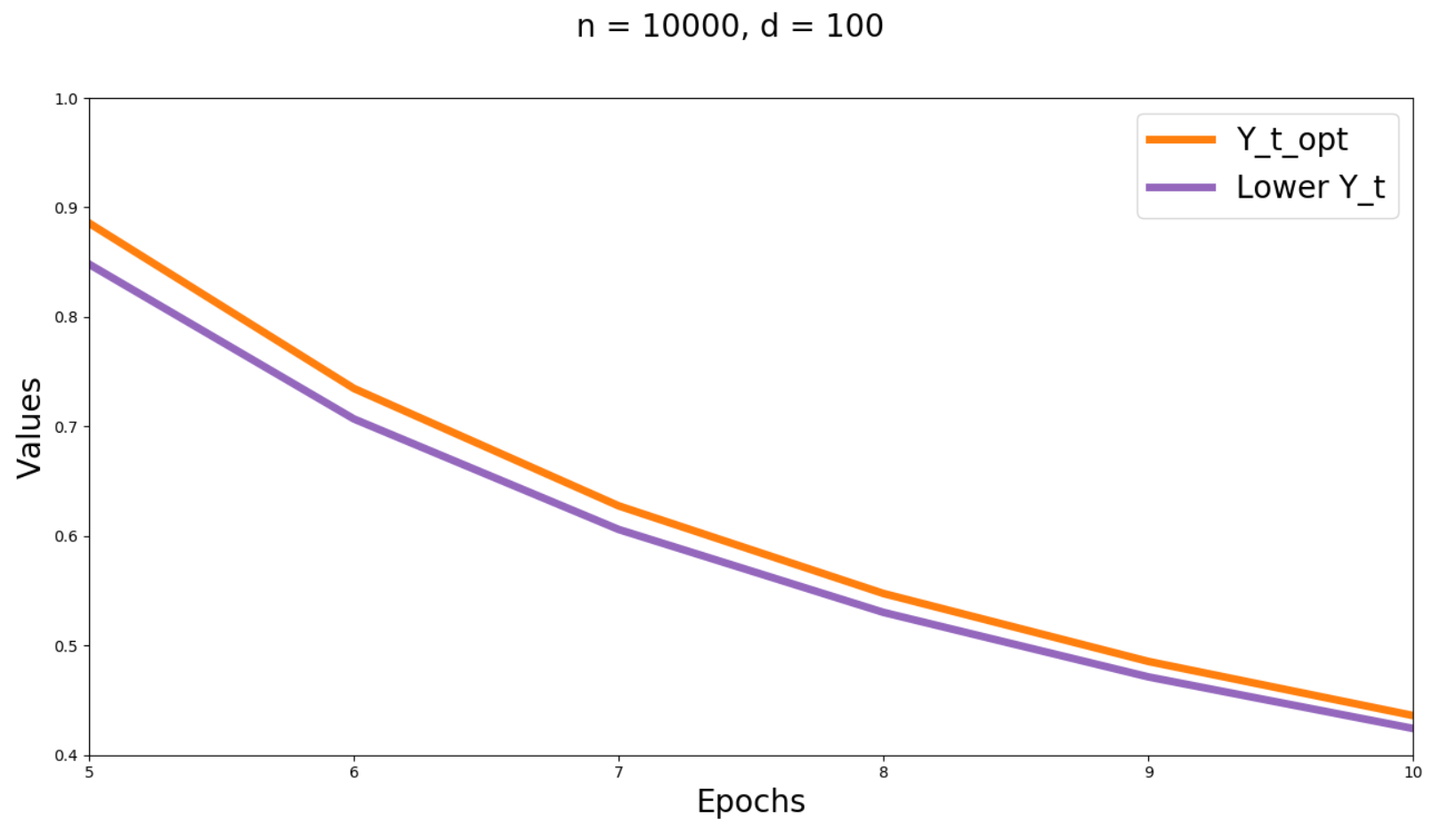}
 \includegraphics[width=0.32\textwidth]{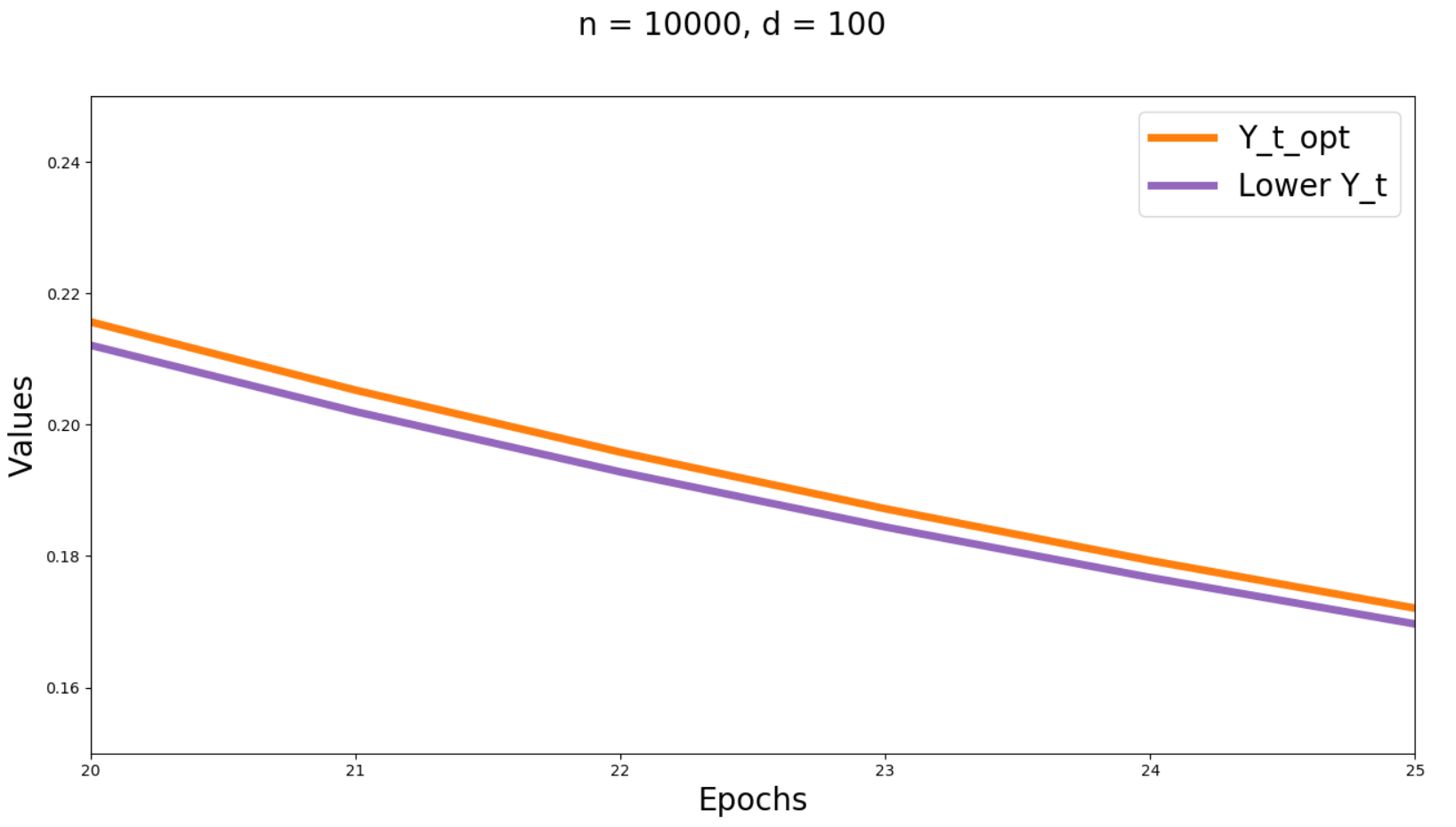} 
   \caption{The difference between ``Lower Y$\_$t'' and ``Y$\_$t$\_$opt'' ($n = 10000$, $d = 100$)}
  \label{xxx}
 \end{figure*}

We denote the labels ``Upper Y$\_$t'' (red line) and ``Lower Y$\_$t'' (violet line) in Figure~\ref{figure_n1} as the upper and lower bounds of $Y_t$ in \eqref{tightUP}  and \eqref{tightLOW} respectively (with a vertical line at epoch $20$ because we expect to see the upper bound take effect when $t \geq T' = 20 L/\mu$,  see supplemental material \ref{supproof}); ``Y$\_$t$\_$opt'' (orange line) as $Y_t$ defined in Theorem~\ref{theorem:example_tightness} computed by using information from  oracle ${\cal U}$; ``Y$\_$t'' (green line) as the squared norm of the difference between $w_t$ and $w_*$, where $w_t$ is generated from Algorithm~\ref{sgd_algorithm} with learning rate  \eqref{optimal_learningrate}. 
Note that $Y_t$ in Figure~\ref{figure_n1} is computed as average of 10 runs of $\|w_t - w_*\|^2$ (not exactly $\mathbb{E}[\|w_t - w_*\|^2]$).  
``Upper Y$\_$t'' (red line), ``Lower Y$\_$t'' (violet line) and 
``Y$\_$t$\_$opt'' (orange line) do not oscillate because they can be correctly computed using formulas \eqref{tightUP}, \eqref{tightLOW}, and \eqref{ytoptimal}, respectively, i.e., they have no variation. The green line ``Y$\_$t'' for stepsize $\eta_t = \frac{2}{\mu t + 4L}$ in Figure~\ref{figure_n1} oscillates because  our analysis does not consider the variance of $\|w_t-w_*\|^2$. 
From (\ref{main_ineq_sgd_new02})
we infer that a decrease in $\eta_t$ leads to a decrease of the variance of $\|w_t-w_*\|^2$.  This fact is reflected in all subfigures in Figure~\ref{figure_n1}. We expect that increasing $d$ and $n$ (the number of dimensions in data and the number of data points) will increase the variance. Hence, it requires larger $t$ to make the variance approach $0$ as shown in Figure~\ref{figure_n1}. For  sufficiently large $t$, the optimality of $\eta_t = \frac{2}{\mu t + 4L}$ is clearly shown in Figure~\ref{figure_n1} when $n=1000$ and $d=10$, i.e., the green line is in between red line (upper bound) and violet line (lower bound). 
We note that ``Lower Y$\_$t'' and ``Y$\_$t$\_$opt'' are very close to each other in Figure~\ref{figure_n1} and the difference between them is shown in Figure~\ref{xxx}. 

\section{Related Work}
\label{sec:related_work}

In~\cite{AgarwalBartlettRavikumarEtAl}, the authors showed that the lower bound of $Y_t$ is $\Ocal(1/t)$ \textit{with} bounded gradient assumption for objective function $F$ over a convex set $\mathcal{S}$. To show the lower bound, the authors use the following three assumptions for the objective function $F$:
\begin{enumerate}
\item  The assumption of a strongly convex objective function, i.e., Assumption~\ref{ass_stronglyconvex} (see Definition 3 in~\cite{AgarwalBartlettRavikumarEtAl}).
\item There exists a bounded convex set $\mathcal{S} \subset \mathbb{R}^d$ such that 
\begin{equation*}
\mathbb{E}[\|\nabla f(w;\xi) \|^2]\leq \sigma^2
\end{equation*}
for all $w\in \mathcal{S} \subset \mathbb{R}^d$ (see Definition 1 in~\cite{AgarwalBartlettRavikumarEtAl}). Notice that this is not the same as the bounded gradient assumption where $\mathcal{S} = \mathbb{R}^d$ is unbounded.
\item The objective function $F$ is a convex Lipschitz function, i.e., there exists a positive number $K$ such that
\begin{equation*}
\|  F(w) - F(w') \| \leq K \|w-w' \|, \forall w, w' \in \mathcal{S} \subset \mathbb{R}^d.
\end{equation*} We notice that this assumption actually implies the assumption on bounded gradients as stated above. 
\end{enumerate}

%
%
%
%

\vspace{2mm}
\noindent 
{\bf On the existence of the assumption of bounded convex set $\mathcal{S} \subset \mathbb{R}^d$ where SGD converges:} let us restate the example in~\cite{Nguyen2018}, i.e. $F(w)=\frac{1}{2}(f_1(w)+f_2(w))$ where $f_1(w)=\frac{1}{2}w^2$ and $f_2(w)=w$. It is obvious that $F$ is strongly convex but $f_1$ and $f_2$ are not. Let $w_0=0 \in \mathcal{S}$, for any number $t\geq 0$, with probability $\frac{1}{2^t}$, the steps of SGD algorithm for all $i<t$ are $w_{i+1}=w_i-\eta_i$. This implies that $w_t=-\sum_{i=1}\eta_i$. Since $\sum_{i=1}\eta_i=\infty$, $w_t$ will escape the set $\mathcal{S}$ when $t$ is sufficiently large. 
We conclude that in ${\cal F}_{str}$ there are objective functions that can escape any bounded set $\mathcal{S}$ with non-zero probability.


If $\mathcal{S}$ is $\mathbb{R}^d$, we have the following results:
\vspace{2mm}

\noindent
{\bf On the non-coexistence of the assumption of a bounded gradient over $\mathbb{R}^d$ and assumption of having strong convexity:}
As pointed out in~\cite{Nguyen2018}, the assumption of bounded gradient does not co-exist with strongly convex assumption. 
As shown in \cite{nesterov2004,bottou2016optimization}, Assumption \ref{ass_stronglyconvex} on strong convexity implies
\begin{align*}
2\mu [ F(w) - F(w_{*}) ] \leq \| \nabla F(w) \|^2 \ , \ \forall w \in \mathbb{R}^d. \tagthis\label{eq:stronglyconvex}
\end{align*}
As shown in~\cite{Nguyen2018}, for any $w\in \mathbb{R}^d$, we have 
\begin{align*}
2\mu[F(w)-F(w_*)] &\overset{(\ref{eq:stronglyconvex})}{\leq}\|\nabla F(w)\|^2 = \|\mathbb{E}[\nabla f(w;\xi)] \|^2 \\
&\leq \mathbb{E}[\|\nabla  f(w;\xi) \|^2] \leq \sigma^2. 
\end{align*}
Therefore, $$ F(w) \leq \frac{\sigma^2}{2\mu} + F(w_*), \forall w \in \mathbb{R}^d.$$
Note that, the from Assumption~\ref{ass_stronglyconvex} and $\nabla F(w_*)=0$, we have
$$
 F(w) \geq \mu \|w-w_*\|^2 + F(w_*), \forall w \in \mathbb{R}^d. 
$$
Clearly, the two last inequalities contradict to each other for sufficiently  large $\|w-w_*\|^2$. Precisely, only when $\sigma$ is equal to $\infty$, then the assumption of bounded gradient and the assumption of strongly convexity of $F$ can co-exist. However, $\sigma$ cannot be $\infty$ and this result implies that there does not exist any objective function $F$ satisfies the assumption of bounded gradients over $\mathbb{R}^d$ and the assumption of having a strongly convex objective function at the same time. 

\vspace{2mm}
\noindent
{\bf On the non-coexistence of the assumption of being convex Lipschitz over $\mathbb{R}^d$ and assumption of being strongly convex:}
Moreover, we can also show that the assumption of convex Lipschitz function does not co-exist with the assumption of being strongly convex. As shown in Section 2.3 in \cite{AgarwalBartlettRavikumarEtAl}, the assumption of Lipschitz function implies that $\|\nabla F(w)\| \leq K, \forall w\in R^d$. Hence, by using the same argument from the analysis of the non-coexistence of bounded gradient assumption and assumption of strongly convex, we can conclude that  these two assumptions cannot co-exist. In other words, there does not exist an objective function $F$ which satisfies the assumption of convex Lipschitz function and assumption of being strongly convex at the same time. 

\subsection{Discussion on the usage of Assumptions in~\cite{AgarwalBartlettRavikumarEtAl}}

As stated in Section 3 and Section 4.1.1 in~\cite{AgarwalBartlettRavikumarEtAl}, the authors construct a class of strongly convex Lipschitz objective function $F$ which has $K=\sigma$. The authors showed that the problem of convex optimization for the constructed class of objective functions $F$ is at least as hard as estimating the biases of $d$ independent coins (i.e., the problem of estimating parameters of Bernoulli variables). As one additional important assumption to prove the lower bound of a first order stochastic algorithm, the authors assume the \textbf{existence} of stepsizes $\eta_t$ which make an first order stochastic algorithm converge for a given objective function $F$ under the three aforementioned assumptions (see Lemma 2 in~\cite{AgarwalBartlettRavikumarEtAl}). Note that the proof of the lower bound of $Y_t$ is described in Theorem 2 in~\cite{AgarwalBartlettRavikumarEtAl} and Theorem 2 uses their Lemma 2. If their Lemma 2 is not valid, then the proof of the lower bound of $Y_t$ in Theorem 2 is also not valid. 

Given the proof strategy in \cite{AgarwalBartlettRavikumarEtAl} of the convergence of a first order stochastic algorithm, one may require that the convex set $\mathcal{S}$ where $F$ has all these nice properties  must be $\mathbb{R}^d$ as explained above. This, however, will lead to the non-coexistence of bounded gradient assumption and strongly convex assumption and the non-coexistence of Lipschitz function assumption and strongly convex assumption as discussed above. In this case, their Lemma 2 is not valid because of non-existence of an objective function $F$, in which case the proof of lower bound of $Y_t$ in Theorem 2 is not correct.

However, we explain why the setup as proposed in~\cite{AgarwalBartlettRavikumarEtAl} may still be acceptable and lead to a proper lower bound: The paper assumes that we only restrict the analysis of SGD in a bounded convex set $\mathcal{S}$ which is not $\mathbb{R}^d$, and only in this bounded set $\mathcal{S}$ we assume that objective function acts like a Lipschitz function (implying bounded gradients in ${\cal S}$). 

There are two possible cases at the $t$-th iteration a first order stochastic algorithm, the algorithm diverges or converges. Let us define $p_t=Pr(w_t \notin \mathcal{S})$. Hence, $Pr(w_t \in \mathcal{S})=1-p_t$. Let $$Y^{conv}_t=\mathbb{E}[\|w_t - w_*\|^2 |w_t \in \mathcal{S}]$$ and $$Y^{div}_t=\mathbb{E}[\|w_t - w_*\|^2 |w_t \notin \mathcal{S}].$$
Since $Y_t =\mathbb{E}[\|w_t - w_*\|^2$, $Y_t$ is equal to 
\begin{align*}
Y_t &= p\cdot Y^{div}_t + (1-p)\cdot Y^{conv}_t\\
				&\geq p\cdot Y^{conv}_t + (1-p)\cdot Y^{conv}_t\\
				&\geq Y^{conv}_t \\
				&\geq \text{ lower bound in~\cite{AgarwalBartlettRavikumarEtAl}}.
\end{align*}

The above derivation hinges on the
 first inequality where we assume $Y^{div}_t\geq Y^{conv}_t$. Typically, for strongly convex objective functions and $w_*,w_0\in \mathcal{S}$), it seems always true that $Y^{div}_t \geq Y^{conv}_t$ because $w_t$ gets far from $w_*$ for the divergence case and it gets close to $w_*$ for the convergence case. Of course a proper proof of this property is still needed in order to rigorously complete the argument leading to the lower bound in \cite{AgarwalBartlettRavikumarEtAl}. In fact this remains an open problem (one can invent strange corner cases that need extra thought/proof).
 
 The above result is interesting because now we \textbf{only} need to prove the convergence of a first order stochastic algorithm in a certain convex set $\mathcal{S}$ with a certain probability $p$. This is completely different from the proof of convergence of e.g. SGD in the general case as in~\cite{Moulines2011} and~\cite{Nguyen2018, gower2019sgd} where we need to prove it with probability 1.

\subsection{Setup}

We describe the setup of the class of strong convex functions proposed in~\cite{AgarwalBartlettRavikumarEtAl}.

As shown in Section 4.1.1~\cite{AgarwalBartlettRavikumarEtAl}, the following two sets are required. 

\begin{enumerate}
\item Subset $\mathcal{V} \subset \{-1,+1 \}^d$ and $\mathcal{V}=\{\alpha^1,\dotsc,\alpha
^M \}$ with $\Delta_H(\alpha^j,\alpha^k) \geq \frac{d}{4}$ for all $j\neq k$, where $\Delta_H$ denotes the Hamming metric, i.e $\Delta_H(\alpha,\beta):=\sum_{i=1}^d \mathbb{I}[\alpha_i \neq \beta_i]$. As discussed by the authors, $|\mathcal{V}|=M\geq (2/\sqrt{e})^{\frac{d}{2}}$. 
\item Subset $\mathcal{F}_{base}=\{f^{+}_i,f^{-}_i,i=1,\dotsc,d\}$ where $f^{+}_i,f^{-}_i$ will be designed depending on the problem at hand. 
\end{enumerate}

Given $\mathcal{V}$, $\mathcal{F}_{base}$ and a constant $\delta \in (0,\frac{1}{4}]$, we define the function class $\mathcal{F}(\delta):=\{F_\alpha, \alpha \in \mathcal{V}\}$ where  
\begin{equation}
\label{eq:Falpha}
F_\alpha(w): = \frac{c}{d}\sum_{i=1}^d \{(1/2+\alpha_i\delta)f^{+}_i(w) + (1/2-\alpha_i\delta)f^{-}_i(w)\}.
\end{equation}
The $\mathcal{F}_{base}$ and constant $c$ are chosen in such a way that $\mathcal{F}(\delta) \subset \mathcal{F}$ where $\mathcal{F}$ is the class of strongly convex objective functions defined over set $\mathcal{S}$ and satisfies all the assumptions as mentioned before. In case $\mathcal{F}$ is the class of strongly convex functions, the key idea to compute the lower bound of SGD proposed in~\cite{AgarwalBartlettRavikumarEtAl} by applying Fano's inequality~\cite{Yu1997} and Le Cam's bound \cite{CoverThomas1991,LeCamothers1973} is as follows: If an SGD algorithm  $\mathcal{M}_t$ works well for optimizing a given function $F_{\alpha^*}, \alpha^* \in \mathcal{V}$ with a given oracle $\mathcal{U}$, then there exists a hypothesis test finding $\hat{\alpha}$ such that:
\begin{equation}
\label{eq:bound_eq}
  \frac{1}{3} \geq \mathrm{Pr}_{\mathcal{U}}[\hat{\alpha}(\mathcal{M}_t)\neq \alpha] \geq 1 -2\frac{16dt\delta^2+\log(2)}{d\log (2/\sqrt{e})}. 
\end{equation}

From (\ref{eq:bound_eq}), we have
\begin{equation*}
\frac{16dt\delta^2+\log(2)}{d\log (2/\sqrt{e})} \approx \frac{16dt\delta^2}{d\log (2/\sqrt{e})} \geq 2/3.
\end{equation*}

Hence, 
\begin{equation}
\label{eq:equa_xxxx}
t \geq \frac{\log(2/\sqrt{e})}{48} \frac{1}{\delta^2}. 
\end{equation}

As shown in Section 4.3~\cite{AgarwalBartlettRavikumarEtAl}, to proceed the proof, we set $Y_t =\frac{c\delta^2r^2}{18(1-\theta)}$. 
Combining with (\ref{eq:equa_xxxx}) yields 
\begin{equation}
\label{eq:equa_xxy}
Y_t \geq \frac{1}{t}\frac{\log(2/\sqrt{e})}{864}\frac{cr^2}{1-\theta}.
\end{equation} 

In addition to the proof of the lower bound, we also need to set $c=\frac{Ld}{rd^{1/p}}$ and $\mu^2 = \frac{L}{rd^{1/p}}(1-\theta)$ where $\mathcal{S}=\mathrm{B}_\infty(r)$. By substituting $c$ and $\mu^2$ into (\ref{eq:equa_xxy}), we obtain:
\begin{equation}
\label{eq:equa_xxzt}
Y_t \geq \frac{1}{t}\frac{\log(2/\sqrt{e})}{864d}\frac{1}{\mu^2}c^2r^2. 
\end{equation} 

To complete the description of the setup in~\cite{AgarwalBartlettRavikumarEtAl}, we briefly describe the proposed oracle $\mathcal{U}$ which outputs some information to the SGD algorithm at each iteration for constructing the stepsize $\eta_t$. There are two types of oracle $\mathcal{U}$ defined as follows.
\begin{enumerate}
\item Oracle $\mathcal{U}_A$: 1-dimensional unbiased gradients
  \begin{enumerate}
   \item Pick an index $i\in {1,\dotsc,d}$ uniformly at random.
   \item Draw $b_i \in \{0,1\}$ according to a Bernoulli distribution with parameter $1/2+\alpha_i\delta$.
    \item For the given input $x\in \mathcal{S}$, return the value $f_i$ and subgradient $\nabla f_i$ of the function $$f_{i,A}:= c[b_i f^{+}_i + (1-b_i)f^{-}_i].$$
  \end{enumerate}
\item Oracle $\mathcal{U}_B$: $d$-dimensional unbiased gradients. 
   \begin{itemize}
    \item For $i=1,\dotsc,d$, draw $b_i \in \{0,1\}$ according to a Bernoulli distribution with parameter $1/2+\alpha_i \delta$.
     \item For the given input $x\in \mathcal{S}$, return the value $f_i$ and subgradient $\nabla f_i$ of the function $$f_{i,B}:= \frac{c}{d} \sum_{i=1}^d [b_i f^{+}_i + (1-b_i)f^{-}_i].$$
   \end{itemize}     
\end{enumerate}

\subsection{Analysis and Comparison}
In this section, we want to compare our lower bound  ($\approx \frac{N}{2\mu^2 t}$) with the one in (\ref{eq:equa_xxzt}) when $t$ is sufficiently large. In order to do this, we need to compute $N=2\mathbb{E}[\|\nabla f(w^*;\xi)\|^2]$ for the strongly convex function class proposed in~\cite{AgarwalBartlettRavikumarEtAl}. For the strongly convex case, the authors defined the base functions as follows. Given a parameter $\theta \in [0,1)$, we have
\begin{align*}
f^{+}_i(w) &= r\theta|w_i + r | + \frac{1-\theta}{4}(w_i+r)^2, \\
f^{-}_i(w) &= r\theta|w_i - r | + \frac{1-\theta}{4}(w_i-r)^2,
\end{align*}
 where $w=(w_1,\dotsc,w_d)$. Let $e_i$ be $1/2+\alpha_i\delta$. Substituting $e_i$ in (\ref{eq:Falpha}) yields $F_\alpha(w)=\frac{1}{d}[\sum_{i=1}^d f_{\alpha,i}(w)]$ where $f_{\alpha,i}(w) = c[e_i f^{+}_i(w) + (1-e_i)f^{-}_i(w)]$. Due to the construction of $F_\alpha$, the definition of $f_{\alpha,i}(w)$ and the construction of oracle $\mathcal{U}_A$ or oracle $\mathcal{U}_B$,  $w^*$ of $F_\alpha$ can be found by finding each $w^*_i$ for each $f_{\alpha,i}(w)$ first. Precisely, we have the following cases:
 
\begin{enumerate}
\item $w_i<-r$: we have 
	     \begin{itemize}
	        \item $f_{\alpha,i}(w)= -r\theta(w_i+r)e_i + \frac{1-\theta}{4}(w_i+r)^2e_i - r\theta (w_i -r)(1-e_i) + \frac{1-\theta}{4}(w_i-r)^2(1-e_i)$.
	        \item $\nabla f_{\alpha,i}(w)= (1-\theta)e_ir -\frac{1+\theta}{2}r + \frac{1-\theta}{2} w_i$.    
	        \item $\nabla  f_{\alpha,i}(w)=0$ at $w^{-r}_i=r[1-2e_i+\frac{2\theta}{1-\theta}]$.
	     \end{itemize}
\item $-r\leq w_i\leq r$: we have 
	     \begin{itemize}
	        \item $f_{\alpha,i}(w)= r\theta(w_i+r)e_i + \frac{1-\theta}{4}(w_i+r)^2e_i - r\theta (w_i -r)(1-e_i) + \frac{1-\theta}{4}(w_i-r)^2(1-e_i)$.
	        \item $\nabla f_{\alpha,i}(w)= (1+\theta)e_ir -\frac{1+\theta}{2}r + \frac{1-\theta}{2} w_i$.    
	        \item $\nabla  f_{\alpha,i}(w)=0$ at $w^{[-r,r]}_i=r\frac{1+\theta}{1-\theta}(1-2e_i)$.
	     \end{itemize}	     
\item $r\leq w_i\leq \infty$: we have 
	     \begin{itemize}
	        \item $f_{\alpha,i}(w)=r\theta(w_i+r) e_i + \frac{1-\theta}{4}(w_i+r)^2e_i + r\theta (w_i -r)(1-e_i) + \frac{1-\theta}{4}(w_i-r)^2(1-e_i)$.
	        \item $\nabla f_{\alpha,i}(w)= (1-\theta)e_ir +\frac{3\theta-1}{2}r + \frac{1-\theta}{2} w_i$.    
	        \item $\nabla  f_{\alpha,i}(w)=0$ at $w^r_i=r[1-2e_i - 2\frac{\theta}{1-\theta}]$.
	     \end{itemize}	     
\end{enumerate}  

Now, we have five important points $w^{-r}_i, w^{[-r,r]}_i, w^r_i, -r$ and $r$ and at these points $F_{\alpha}$ can be minimum. We consider the following cases

\begin{enumerate}
\item $\alpha_i=-1$ and then $e_i=\frac{1}{2}+\alpha_i \delta=\frac{1}{2}-\delta$ where $\delta \in [0,1/4)$, we have 
     \begin{itemize}
       \item $w^{-r}_i=r[\frac{2\theta}{1-\theta}+2\delta]>-r$.
       \item $w^{[-r,r]}_i= r\frac{1+\theta}{1-\theta}(2\delta)$. In this case $w^{[-r,r]}_i$ may belong $[-r,r]$ or it may be greater than $r$.
       \item $w^r_i=r(2\delta - \frac{2\theta}{1-\theta})<r$ . 
     \end{itemize}
     This result implies $F_{\alpha}$ is minimum at $w^*_i=r$ and $\nabla f_{\alpha,i}(w^*)=cr[(1-\theta)e_i + \theta]=cr[(1-\theta)(1/2-\delta)+ \theta]$. Or it can be minimum at $w^{[-r,r]}_i$ if $w^{[-r,r]}_i \in [-r,r]$ and $\nabla f_{\alpha,i}(w^*)=0$.
\item $\alpha_i=+1$ and then $e_i=\frac{1}{2}+\alpha_i \delta=\frac{1}{2}+\delta$ where $\delta \in [0,1/4)$, we have 
     \begin{itemize}
       \item $w^{-r}_i=r[\frac{2\theta}{1-\theta}-2\delta]$. Since $\frac{2\theta}{1-\theta}-2\delta>-1$ when $\delta \in [0,1/4)$ and $\theta \in [0,1)$. Hence $w^{-r}_i>-r$.
       \item $w^{[-r,r]}_i= r\frac{1+\theta}{1-\theta}(-2\delta)<0$. In this case $w^{[-r,r]}_i$ may belong $[-r,r]$ or it may be smaller than $-r$.
       \item $w^r_i=r(-2\delta - \frac{2\theta}{1-\theta})<r$. 
     \end{itemize}
     This result implies $F_{\alpha}$ is minimum at $w^*_i=-r$ and $\nabla f_{\alpha,i}(w^*)=cr[(1-\theta)e_i -1]=cr[(1-\theta)(1/2+\delta)-1]$. Or it can be minimum at $w^{[-r,r]}_i$ if $w^{[-r,r]}_i \in [-r,r]$ and $\nabla f_{\alpha,i}(w^*)=0$.    
\end{enumerate} 
 
By definition, we have 

\begin{align*}
N&=2\mathbb{E}[\|\nabla f_i(w^*) \|^2]= 2\frac{1}{d}\sum_{i=1}^d [e_i\|c\nabla f^{+}_i(w^*) \|^2 +(1-e_i)\|c\nabla f^{-}_i(w^*) \|^2]
\end{align*}

From the analysis above, we have four possible $w^*_i$, i.e., $-r,r,r\frac{1+\theta}{1-\theta}(-2\delta)$ and $r\frac{1+\theta}{1-\theta}(2\delta)$.  If we plug $w^*$ which has $w^*_i=-r$ or $w^*_i=r$, then we have $[e_i\|c\nabla f^{+}_i(w^*) \|^2 +(1-e_i)\|c\nabla f^{-}_i(w^*) \|^2]=(1/2 -\delta)c^2r^2$. For $w^*_i$ which has $w^*_i=r\frac{1+\theta}{1-\theta}(-2\delta)$ or $r\frac{1+\theta}{1-\theta}(2\delta)$,  we have   $[e_i\|c\nabla f^{+}_i(w^*) \|^2 +(1-e_i)\|c\nabla f^{-}_i(w^*) \|^2]=(1/4 -\delta^2)(1+\theta)^2c^2r^2$. This proves that
$$N=2 \beta c^2 r^2$$
with $\beta$ somewhere in the range
$$ [ (\frac{1}{2} - \delta), (\frac{1}{4}-\delta^2) (1+\theta)^2 ] \mbox{ or } 
[(\frac{1}{4}-\delta^2) (1+\theta)^2 , (\frac{1}{2} - \delta)],$$
where $\delta\in [0,1/4)$ and $\theta \in [0,1)$.

Substituting $N=2\beta c^2r^2$ into (\ref{eq:equa_xxzt}) yields 
\begin{equation}
\label{eq:equa_xxzww}
Y_t \geq \frac{\log(2/\sqrt{e})}{(864 \cdot d)(2\beta)}\frac{N}{\mu^2 t}, 
\end{equation} 
which is further minimized by taking
$$ \beta = \max \{ (\frac{1}{2} - \delta), (\frac{1}{4}-\delta^2) (1+\theta)^2\} .$$
Notice that, given our freedom in choosing $\delta$ and $\theta$, we can minimize $\beta$ as a function of $\delta$ and $\theta$ in order to maximize the lower bound in (\ref{eq:equa_xxzww}). This gives (in the limit) $\delta=1/4$ with $\theta\leq 2/\sqrt{3} -1= 0.155$ leading to $\beta=1/4$. This leads to the final lower bound
$$
Y_t \geq \frac{\log(2/\sqrt{e})}{432 \cdot d}\frac{N}{\mu^2 t}.
$$

Clearly, the lower bound in 
is much smaller than our lower bound of $\approx \frac{N}{2\mu^2 t}$ when $t$ is sufficiently large. Moreover, this lower bound depends on $1/d$ and it becomes smaller when $d$ increases. 

\end{document}